\documentclass[11pt]{article}
\usepackage{color}
\usepackage{amsfonts,amssymb,a4, graphicx, bm, makeidx}
\usepackage{amsmath}
\usepackage{soul}
\usepackage{mathrsfs}
\usepackage[all]{xy}
\usepackage[latin1]{inputenc}
\usepackage[dvips]{geometry}
\def\theequation{\thesection.\arabic{equation}}

\newcommand{\qed}{\hfill\rule{3mm}{3mm}}
\newtheorem{corollary}{Corollary}[section]
\newtheorem{theorem}[corollary]{Theorem}
\newtheorem{lemma}[corollary]{Lemma}
\newtheorem{proposition}[corollary]{Proposition}
\newtheorem{definition}[corollary]{Definition}

\makeatletter \@addtoreset{equation}{section} \makeatother
\definecolor{GreenYellow}{cmyk}{0.15,0,0.69,0}
\definecolor{Yellow}{cmyk}{0,0,1,0}
\definecolor{Goldenrod}{cmyk}{0,0.10,0.84,0}
\definecolor{Dandelion}{cmyk}{0,0.29,0.84,0}
\definecolor{Apricot}{cmyk}{0,0.32,0.52,0}
\definecolor{Peach}{cmyk}{0,0.50,0.70,0}
\definecolor{Melon}{cmyk}{0,0.46,0.50,0}
\definecolor{YellowOrange}{cmyk}{0,0.42,1,0}
\definecolor{Orange}{cmyk}{0,0.61,0.87,0}
\definecolor{BurntOrange}{cmyk}{0,0.51,1,0}
\definecolor{Bittersweet}{cmyk}{0,0.75,1,0.24}
\definecolor{RedOrange}{cmyk}{0,0.77,0.87,0}
\definecolor{Mahogany}{cmyk}{0,0.85,0.87,0.35}
\definecolor{Maroon}{cmyk}{0,0.87,0.68,0.32}
\definecolor{BrickRed}{cmyk}{0,0.89,0.94,0.28}
\definecolor{Red}{cmyk}{0,1,1,0}
\definecolor{OrangeRed}{cmyk}{0,1,0.50,0}
\definecolor{RubineRed}{cmyk}{0,1,0.13,0}
\definecolor{WildStrawberry}{cmyk}{0,0.96,0.39,0}
\definecolor{Salmon}{cmyk}{0,0.53,0.38,0}
\definecolor{CarnationPink}{cmyk}{0,0.63,0,0}
\definecolor{Magenta}{cmyk}{0,1,0,0}
\definecolor{VioletRed}{cmyk}{0,0.81,0,0}
\definecolor{Rhodamine}{cmyk}{0,0.82,0,0}
\definecolor{Mulberry}{cmyk}{0.34,0.90,0,0.02}
\definecolor{RedViolet}{cmyk}{0.07,0.90,0,0.34}
\definecolor{Fuchsia}{cmyk}{0.47,0.91,0,0.08}
\definecolor{Lavender}{cmyk}{0,0.48,0,0}
\definecolor{Thistle}{cmyk}{0.12,0.59,0,0}
\definecolor{Orchid}{cmyk}{0.32,0.64,0,0}
\definecolor{DarkOrchid}{cmyk}{0.40,0.80,0.20,0}
\definecolor{Purple}{cmyk}{0.45,0.86,0,0}
\definecolor{Plum}{cmyk}{0.50,1,0,0}
\definecolor{Violet}{cmyk}{0.79,0.88,0,0}
\definecolor{RoyalPurple}{cmyk}{0.75,0.90,0,0}
\definecolor{BlueViolet}{cmyk}{0.86,0.91,0,0.04}
\definecolor{Periwinkle}{cmyk}{0.57,0.55,0,0}
\definecolor{CadetBlue}{cmyk}{0.62,0.57,0.23,0}
\definecolor{CornflowerBlue}{cmyk}{0.65,0.13,0,0}
\definecolor{MidnightBlue}{cmyk}{0.98,0.13,0,0.43}
\definecolor{NavyBlue}{cmyk}{0.94,0.54,0,0}
\definecolor{RoyalBlue}{cmyk}{1,0.50,0,0}
\definecolor{Blue}{cmyk}{1,1,0,0}
\definecolor{Cerulean}{cmyk}{0.94,0.11,0,0}
\definecolor{Cyan}{cmyk}{1,0,0,0}
\definecolor{ProcessBlue}{cmyk}{0.96,0,0,0}
\definecolor{SkyBlue}{cmyk}{0.62,0,0.12,0}
\definecolor{Turquoise}{cmyk}{0.85,0,0.20,0}
\definecolor{TealBlue}{cmyk}{0.86,0,0.34,0.02}
\definecolor{Aquamarine}{cmyk}{0.82,0,0.30,0}
\definecolor{BlueGreen}{cmyk}{0.85,0,0.33,0}
\definecolor{Emerald}{cmyk}{1,0,0.50,0}
\definecolor{JungleGreen}{cmyk}{0.99,0,0.52,0}
\definecolor{SeaGreen}{cmyk}{0.69,0,0.50,0}
\definecolor{Green}{cmyk}{1,0,1,0}
\definecolor{ForestGreen}{cmyk}{0.91,0,0.88,0.12}
\definecolor{PineGreen}{cmyk}{0.92,0,0.59,0.25}
\definecolor{LimeGreen}{cmyk}{0.50,0,1,0}
\definecolor{YellowGreen}{cmyk}{0.44,0,0.74,0}
\definecolor{SpringGreen}{cmyk}{0.26,0,0.76,0}
\definecolor{OliveGreen}{cmyk}{0.64,0,0.95,0.40}
\definecolor{RawSienna}{cmyk}{0,0.72,1,0.45}
\definecolor{Sepia}{cmyk}{0,0.83,1,0.70}
\definecolor{Brown}{cmyk}{0,0.81,1,0.60}
\definecolor{Tan}{cmyk}{0.14,0.42,0.56,0}
\definecolor{Gray}{cmyk}{0,0,0,0.50}
\definecolor{Black}{cmyk}{0,0,0,1}
\definecolor{White}{cmyk}{0,0,0,0}

\begin{document}
\def\theequation{\thesection.\arabic{equation}}

\def\blu{\color{Blue}}
\def\mag{\color{Maroon}}
\def\red{\color{Red}}
\def\green{\color{ForestGreen}}
\def\prob{{\rm Prob}}

\def\reff#1{(\protect\ref{#1})}

\let\a=\alpha \let\b=\beta \let\ch=\chi \let\d=\delta \let\e=\varepsilon
\let\f=\varphi \let\g=\gamma \let\h=\eta    \let\k=\kappa \let\l=\lambda
\let\m=\mu \let\n=\nu \let\o=\omega    \let\p=\pi \let\ph=\varphi
\let\r=\rho \let\s=\sigma \let\t=\tau \let\th=\vartheta
\let\y=\upsilon \let\x=\xi \let\z=\zeta
\let\D=\Delta \let\F=\Phi \let\G=\Gamma \let\L=\Lambda \let\Th=\Theta
\let\O=\Omega \let\P=\Pi \let\Ps=\Psi \let\Si=\Sigma \let\X=\Xi
\let\Y=\Upsilon

\global\newcount\numsec\global\newcount\numfor
\gdef\profonditastruttura{\dp\strutbox}
\def\senondefinito#1{\expandafter\ifx\csname#1\endcsname\relax}
\def\SIA #1,#2,#3 {\senondefinito{#1#2}
\expandafter\xdef\csname #1#2\endcsname{#3} \else
\write16{???? il simbolo #2 e' gia' stato definito !!!!} \fi}
\def\etichetta(#1){(\veroparagrafo.\veraformula)
\SIA e,#1,(\veroparagrafo.\veraformula)
 \global\advance\numfor by 1
 \write16{ EQ \equ(#1) ha simbolo #1 }}
\def\etichettaa(#1){(A\veroparagrafo.\veraformula)
 \SIA e,#1,(A\veroparagrafo.\veraformula)
 \global\advance\numfor by 1\write16{ EQ \equ(#1) ha simbolo #1 }}
\def\BOZZA{\def\alato(##1){
 {\vtop to \profonditastruttura{\baselineskip
 \profonditastruttura\vss
 \rlap{\kern-\hsize\kern-1.2truecm{$\scriptstyle##1$}}}}}}
\def\alato(#1){}
\def\veroparagrafo{\number\numsec}\def\veraformula{\number\numfor}
\def\Eq(#1){\eqno{\etichetta(#1)\alato(#1)}}
\def\eq(#1){\etichetta(#1)\alato(#1)}
\def\Eqa(#1){\eqno{\etichettaa(#1)\alato(#1)}}
\def\eqa(#1){\etichettaa(#1)\alato(#1)}
\def\equ(#1){\senondefinito{e#1}$\clubsuit$#1\else\csname e#1\endcsname\fi}
\let\EQ=\Eq

\def\pp{{\bm p}}\def\pt{{\tilde{\bm p}}}


\def\\{\noindent}
\let\io=\infty
\def\ee{\end{equation}}
\def\be{\begin{equation}}

\def\VU{{\mathbb{V}}}
\def\EE{{\mathbb{E}}}
\def\N{\mathbb{N}}
\def\U{\mathbb{U}}
\def\GI{{\mathbb{G}}}
\def\TT{{\mathbb{T}}}
\def\C{\mathbb{C}}
\def\CC{{\mathcal C}}
\def\KK{{\mathcal K}}
\def\II{{\mathcal I}}
\def\LL{{\cal L}}
\def\RR{{\cal R}}
\def\SS{{\cal S}}
\def\NN{{\cal N}}
\def\HH{{\cal H}}
\def\GG{{\cal G}}
\def\PP{{\cal P}}
\def\AA{{\cal A}}
\def\BB{{\cal B}}
\def\FF{{\cal F}}
\def\v{\vskip.1cm}
\def\vv{\vskip.2cm}
\def\gt{{\tilde\g}}
\def\E{{\mathcal E} }
\def\EI{{\mathbb E} }
\def\I{{\rm I}}
\def\rfp{R^{*}}
\def\rd{R^{^{_{\rm D}}}}
\def\ffp{\varphi^{*}}
\def\ffpt{\widetilde\varphi^{*}}
\def\fd{\varphi^{^{_{\rm D}}}}
\def\fdt{\widetilde\varphi^{^{_{\rm D}}}}
\def\pfp{\Pi^{*}}
\def\pd{\Pi^{^{_{\rm D}}}}
\def\pbfp{\Pi^{*}}
\def\fbfp{{\bm\varphi}^{*}}
\def\fbd{{\bm\varphi}^{^{_{\rm D}}}}
\def\rfpt{{\widetilde R}^{*}}
\def\A{{{\mathcal O}}}
\def\ef{\mathfrak{f}}
\def\Ti{\mathfrak{T}}
\def\Mi{\mathfrak{M}}

\def\tende#1{\vtop{\ialign{##\crcr\rightarrowfill\crcr
              \noalign{\kern-1pt\nointerlineskip}
              \hskip3.pt${\scriptstyle #1}$\hskip3.pt\crcr}}}
\def\otto{{\kern-1.truept\leftarrow\kern-5.truept\to\kern-1.truept}}
\def\arm{{}}
\font\bigfnt=cmbx10 scaled\magstep1

\newcommand{\card}[1]{\left|#1\right|}
\newcommand{\und}[1]{\underline{#1}}
\def\1{\rlap{\mbox{\small\rm 1}}\kern.15em 1}
\def\ind#1{\1_{\{#1\}}}
\def\bydef{:=}
\def\defby{=:}
\def\buildd#1#2{\mathrel{\mathop{\kern 0pt#1}\limits_{#2}}}
\def\card#1{\left|#1\right|}
\def\proof{\noindent{\bf Proof. }}
\def\qed{ \square}
\def\trp{\mathbb{T}}
\def\trt{\mathcal{T}}
\def\Z{\mathbb{Z}}
\def\be{\begin{equation}}
\def\ee{\end{equation}}
\def\bea{\begin{eqnarray}}
\def\eea{\end{eqnarray}}
\def\kk{{\bf k}}
\def\Ti{\mathfrak{T}}
\def\Mi{\mathfrak{M}}
\def\begn{\begin{aligned}}
\def\egn{\end{aligned}}
\def\ti{{\rm\bf  t}}\def\mi{{\rm\bf m}}
\def\Va{{V^a_{\rm h.c.}}}
\def\Re{{\mathbb{R}}}
\def\T{{\mathcal{T}}}
\def\hL{{\L}}
\def\ev{\mathfrak{e}}
\def\obj{{\rm supp}}\def\fa{\FF}
\def\E{{\cal E}}
\def\EA{{E_\A}}
\def\0{\emptyset}
\def\Ni{\overline{\N}}

\title{On the zero-free region  for the chromatic polynomial of graphs with maximum degree $\D$ and girth $g$}

\author{
\\
Paula M. S. Fialho$^1$,  Emanuel Juliano$^1$, Aldo Procacci$^2$. \\
\\
\small{$^1$Departamento de Ci\^encia da Computa\c{c}\~ao UFMG, }
\small{30161-970 - Belo Horizonte - MG
Brazil}\\
\small{$^2$Departamento de Matem\'atica UFMG,}
\small{ 30161-970 - Belo Horizonte - MG
Brazil}\\
}

\maketitle

\begin{abstract}
\\The purpose  of the present paper is to provide, for all pairs of integers $(\D,g)$  with  $\D\ge 3$ and $g\ge 3$, a positive number  $C(\D, g)$ such that
chromatic polynomial $P_\GI(q)$ of a graph with maximum degree $\D$ and finite girth $g$ is free of zero if $|q|\ge C(\D, g)$.
Our bounds
enlarge  the zero-free region in the complex plane of  $P_\GI(q)$  in comparison to  previous bounds.
In particular,
for small values  of $\D$  our estimates  yield a sensible improvement  on the bounds  recently obtained by Jenssen, Patel and Regts in \cite{JPR}, while they coincide with those of \cite{JPR} when $\D\to \infty$.

\end{abstract}

\section{Introduction and results}

\\Given a  graph $\GI=(\VU,\EE)$, a coloring  of $\VU$ with $q\in \N$ colors is a function $\k: \VU\to [q]$. A coloring  of $\VU$ is called {\it proper} if for any edge $\{x,y\}\in \EE$ it holds that
$\k(x)\neq\k(y)$. Letting  $\KK^*_\VU(q)$ be the set of all proper colorings of $\VU$ with $q$ colors, the quantity $P_\mathbb{G}(q):= |\KK^*_\VU(q)|$, i.e. the number of proper colorings with $q$ colors of the graph $\mathbb{G}$, is, as a function of $q$, a polynomial known as the \emph{chromatic polynomial} of $\GI$.

\\It is long known (see \cite{KE}, \cite{sok01} and references therein) that  $P_\mathbb{G}(q)$
coincides with the partition function of the antiferromagnetic
Potts model with $q$ states on $\mathbb{G}$ at zero temperature, which in turn can be written in terms of a ``polymer gas" partition function
 (see e.g. Proposition 2.1 in \cite{sok01}). Namely,
\be\label{zqx}
P_\GI(q)=q^{|\VU|}\Xi_\GI(q)
\ee
where $\Xi_\GI(q)$ is the grand canonical partition function of a ``polymer gas'', in which the ``polymers'' are  subsets $R\subset \mathbb{V}$, with cardinality $|R|\ge 2$,  subjected to a non-overlapping constraint and with activities $z(R,q)$ that depend on the topological structure
of $\mathbb{G}$. Explicitly,
\be\label{XiG}
\Xi_\GI(q) =1+\sum_{k\ge 1}\sum_{\{R_1,\dots, R_k\}:\, R_i\subset  \VU\atop |R_i|\ge 2,\;
R_i\cap R_j=\emptyset} \prod_{i =1}^{k}z(R_i,q)
\ee
with
\be\label{actrq}
z(R,q)={1\over q^{|R|-1}} \sum_{g\in\mathcal{C}_{\GI|_R}}(-1)^{|E_g|}
\ee
where $\GI|_R$ is the restriction of $\GI$ to $R$ (i.e. $\GI|_R=(R, \EE|_R)$ where $\EE|_R=\{\{x,y\}\subset R: \{x,y\}\in \EE\}$) and $\mathcal{C}_{\GI|_R}$ represents the set of all connected spanning subgraphs of $\GI|_R$.

\\Therefore to find zero-free regions of $P_\GI(q)$ when $q\in \C$ it is equivalent to find  zero-free regions of the partition
function $\Xi_\GI(q)$ of the above polymer gas. The latter is  a classical problem in statistical mechanics, deeply related to phase transitions.

\\For a general finite graph $\GI$,
whose sole information is its maximum  degree $\D$,  Sokal \cite{sok01} proved that all the zeros of $P_\GI(q)$ lie
in the disc $|q| < C_{\rm So}(\D)$, where $C_{\rm So}(\D)$ is an explicitly computable function of $\D$
(see \cite{sok01} Table 1, some values for $C_{\rm So}(\D)$). He also proved that  $C_{\rm So}(\D)\le 7.963906\D$ for any $\D$. The strategy followed by Sokal to prove his results was to find values of $q\in \C$ such that the absolute value of the logarithm of the partition function $\Xi_\GI(q)$  defined in  \reff{XiG} is bounded (and therefore $\Xi_\GI(q)\neq \0$). To this end, he  used  the so-called Dobrushin criterion \cite{dob96} for the convergence of the cluster expansion of the abstract polymer gas. Later, following  the same strategy developed by Sokal, Fern\'andez and Procacci \cite{FP2} improved the results  of Sokal using their improved version of the Dobrushin criterion given in \cite{FP}, and bounding the activities $z(q,R)$  given in \reff{actrq} using the so called Penrose identity \cite{pen67} (see also ahead). In \cite{FP2} the authors provide a function $ C_{\rm FP}(\D)$ (also explicitly computable)
which is, for any $\D$, smaller than the Sokal function $C_{\rm So}(\D)$ and bounded above by $6.9077\D$. Moreover, $C_{\rm FP}(\D)$
can be done even smaller by using some further information about the structure of the graph $\GI$ beyond its maximum degree.
(see Section 2 of \cite{FP2},   Table 1).

\\Very recently, Jenssen, Patel and Regts \cite{JPR}, using a  representation of $P_\GI(q)$
known as  Whitney broken circuit Theorem \cite{W},
have improved the bound given in \cite{FP2} for the zero-free region of the chromatic polynomial. In \cite{JPR} the  authors show, by an inductive method {\it a la} Dobrushin  inspired by the proof of Proposition 3.1 in \cite{BFP}, that for any
graph $\GI$ with maximum degree $\D$ and girth $g\ge 3$ the chromatic polynomial $P_\GI(q)$
is free of zero   outside the disk  $|q|<K_g \Delta$,
where
\be\label{chat}
K_g=\min_{a\in (0,1)}\left({b_g(a)\over (1-a)\ln b_g(a)}\right)
\ee
with $b_g(a)$ being the unique solution of the equation
\be\label{chat1}
\exp\left\{{(1-a)\ln b\over b}\right\} -1+  {b(\ln b)^{g-1}\over  2(1-{\ln b})}=a
\ee
in the interval $b\in [1,e)$.
By the above definitions, $K_g$ is  a decreasing and explicitly computable function of $g$
(see Figure  1 in \cite{JPR} for some values of $K_g$ for  $3\le g\le 100$).  In particular, via \reff{chat} and \reff{chat1}, one gets   that $K_3= 5.93148$, which improves  the asymptotic bound  $C_{\rm FP}(\D)\le 6.9077\D$  obtained in \cite{FP2}, and  that $\lim_{g\to \infty} K_g= 3.85977$.

\\As
 Jenssen Patel and Regts comment (see Section 5  in \cite{JPR}), the  zero-free region of the chromatic polynomial $P_\GI(q)$ established in their paper  could be improved  for graphs with small maximum degree $\D$.
The purpose of this work is precisely this. Namely, we provide,  for all pairs of integers $(\D,g)$  with $\D\ge 3$ and $g\ge 3$, a positive number  $C(\D, g)$
such that the chromatic polynomial  $P_\GI(q)$ of a graph $\GI$ with maximum degree $\D$ and finite girth $g$ is free of zero if $|q|\geq C(\D, g)$ where the function $C(\D, g)$ is such that $C(\D, g)< K_g\D$  for all $\D$.

\\Our  results  are summarized through the following  theorems.

\begin{theorem}\label{teowhit}
Given a graph  $\GI=(\VU,\EE)$  with  maximum degree $\D\ge 3$ and  finite girth  $g\ge 3$,
there exists an explicitly computable  constant
$C(\D, g)$ such that $P_\GI(q) \not = 0$ for any $q \in \mathbb{C}$ satisfying
\[
|q|\ge C(\D, g).
\]
\end{theorem}

\begin{theorem}\label{coro22}
For any fixed $g\ge 3$,  the sequence $\left\{{C(\D,g)/\D}\right\}_{\D\ge 3}$ is monotonic increasing and
\[
\lim_{\D\to\infty}{C(\D,g)\over \D}\le K_g\label{kstar}
\]
where $K_g$ is the
Jenssen-Patel-Regst
constant defined in \reff{chat}.
\end{theorem}
The analytic procedure that must be followed to explicitly compute  the constant $C(\D,g)$ as a function of $\D$ and $g$ will be given in the next section.

\\The bounds obtained in Theorem \ref{teowhit}
 enlarge  the zero-free region of the chromatic polynomial $P_\GI(q)$  in comparison to the bounds recently given in \cite{JPR}. In particular,
for small values  of $\D$  our estimative  yields a sensible improvement  on \cite{JPR}, while it coincides with that of \cite{JPR} when $\D\to \infty$.

\\The strategy used  in the proof of Theorem \ref{teowhit} follows closely that  developed in \cite{JPR}. In particular,
as in \cite{JPR}, we use the    representation of $\X_\GI(q)$ given by
the   Whitney Broken Circuit Theorem (see ahead) and  then  we revisit and adapt   the inductive strategy developed in \cite{JPR}
to accommodate efficiently the case in which $\D$ is not large.

\vv


\section{The construction of the constant $C(\D,g)$}\label{sec3}

\\Given an integer $\D\ge 3$, let
\be\label{rodel2}
\r_\D= \left(\D-1\over \D-2\right)^\D~~~~{\rm and}~~~~R_\D= {(\D-2)^{\D-2}\over (\D-1)^{\D-1}}~.
\ee
Define, for any $b\in [1,\r_\D)$,
\be\label{hdb}
x_\D(b)=
{b^{2\over \D}-b^{1\over \D}\over b}
\ee
and note that the quantity $x_\D(b)$ is an increasing function of $b\in[1,\r_\D)$, varying from 0 to $R_\D$.

\\Define then,   for any $d\in \N$,   $\D\ge 3$ and  $b\in [1,\r_\D)$,

\be\label{fdelta2}
f^d_\D(b)=
\displaystyle{{1\over (\D-1)}~{[(\D-1)(1-b^{-{1\over \D}})]^{d}\over  [1-(\D-1)(1-b^{-{1\over \D}})]}}~.
 \ee

\\Moreover, given  integers $\D\ge 3$ and $g\ge 3$, we define, for $a \in [0,1]$, $b\in [1, \r_\D)$,  the function
\be\label{acca}
K^g_\D(a,b)=\Big[1+{(1-a)}x_\D(b)\Big]^{\Delta} - 1 +
{f^{g-2}_\D(b)}{\D\choose 2}bx_\D(b)(1+ bx_\D(b))^{\D-1}.
\ee
and we set
\be\label{badg2}
b_\D^g(a)= \sup\{b\in [1,\r_\D): ~K^g_\D(a,b)\le a\},
\ee
\vv
\be\label{zdga}
z_\D^g(a)=(1-a) x_\D(b_\D^g(a)).
\ee
%

\vv

\begin{lemma}\label{lemma:RemarkA}
For any fixed $a\in [0,1]$ and for any $\D\ge 3$, the quantity $b_\D^g(a)$ is  the unique solution of the equation
$K^g_\D(a,b)= a$ as $b$ varies in the interval $[1,\r_\D)$, where $\r_\D$ is defined in (\ref{rodel2}). Moreover, for all $\D\ge 3$ and
all $g\ge 3$ we have that $\max_{a\in [0,1]} {z^g_\D(a)}$ is a well defined number.
\end{lemma}

\\{\bf Proof:}
Note that, for any fixed $a\in [0,1]$ and for any $\D\ge 3$, the function $K^g_\D(a,b)$ given in \reff{acca}
is a continuous increasing  function of  $b\in[1,\r_\D)$ such that $K^g_\D(a,1)=0$ and  $\lim_{b\to \r_\D}K^g_\D(a,b)=+\infty$. Therefore, $b_\D^g(a)$ is the unique solution of the equation $K^g_\D(a,b)= a$ as $b$ varies in the interval $[1,\r_\D)$.
Observe also that function $b_\D^g(a)$ in the domain $a\in [0,1]$ is defined implicitly via the equation $K_\D^g(a,b)-a=0$ with $K_\D^g(a,b)-a$ continuous in $a\times b\in [0,1]\times [1,\r_\D)$. Therefore $b_\D^g(a)$  is continuous in the interval $a\in [0,1]$.
Consequently, $z_\D^g(a)$, as product and composition of continuous functions of $a$, is in turn a continuous function of $a$ in the interval $a\in [0,1]$ and moreover $z_\D^g(0)=z_\D^g(1)=0$. This implies that $\max_{a\in [0,1]} {z^g_\D(a)}$ is a well defined number for all $\D\ge 3$ and all $g\ge 3$.

~~~~~~~~~~~~~~~~~~~~~~~~~~~~~~~~~~~~~~~~~~~~~~~~~~~~~~~~~~~~~~~~~~~~~~~~~~~~~~~~~~~~~~~~~~~~~~~~~~~~~~~~~~~~~~~~~~~~~~~~~~$\Box$
\begin{definition}
Let  $z^g_\D(a)$ be the function of $a\in [0,1]$ defined in \reff{zdga}. Then, for any $\D\ge 3$ and any $g\ge 3$, we define
\be\label{CDg}
C(\D,g)=\left[\max_{a\in [0,1]} {z^g_\D(a)}\right]^{-1}.
\ee
\end{definition}

\\In  Figure 1 we present a  table comparing, for various values of $\D$ and $g$,  the Jenssen-Patel-Regts constant $K_g$
given in \cite{JPR}
with the constant $C(\D,g)/\D$ obtained via \reff{CDg}. In Figure 2 we plot the behaviour of the quantity $C(\D,g)/\D$ as a function of $\D\ge 3$ for  several fixed values of $g$  and  as a function of $g\ge 3$ for  several fixed values of $\D$.
\begin{figure}[!ht]
    \centering
 \resizebox{6.5cm}{!}{   \begin{tabular}{|c|c|c|c|c|c|}
    \hline
        $\Delta$ & $g$ & $a$ & $b^g_\D(a)$ & $C(\Delta, g)/\Delta$ & $K_g$\\ \hline
        3 & 3 & 0.39625 & 1.57848 & 4.55449 & 5.93148 \\ \hline
        4 & 3 & 0.37712 & 1.53409 & 4.89965 & 5.93148 \\ \hline
        5 & 3 & 0.36658 & 1.5107 & 5.10631 & 5.93148 \\ \hline
        6 & 3 & 0.35989 & 1.49623 & 5.24396 & 5.93148 \\ \hline
        20 & 3 & 0.33838 & 1.45155 & 5.72529 & 5.93148 \\ \hline
        3 & 4 & 0.39082 & 1.76326 & 3.83755 & 5.23445 \\ \hline
        3 & 5 & 0.39411 & 1.92992 & 3.48035 & 4.87264 \\ \hline
        3 & 10 & 0.41768 & 2.59496 & 2.88884 & 4.65234\\ \hline
        3 & 25 & 0.44373 & 3.80294 & 2.60289 & 3.97497\\ \hline
        3 & 100 & 0.45824 & 5.81488 & 2.49247 & 3.87487\\ \hline
    \end{tabular}
    }
    \caption{\footnotesize{A comparison between bounds of ref. \cite{JPR} and those of Theorem \ref{teowhit}}}
\end{figure}
\begin{figure}[!ht]
\begin{center}
\includegraphics[width=9cm,height=6cm]{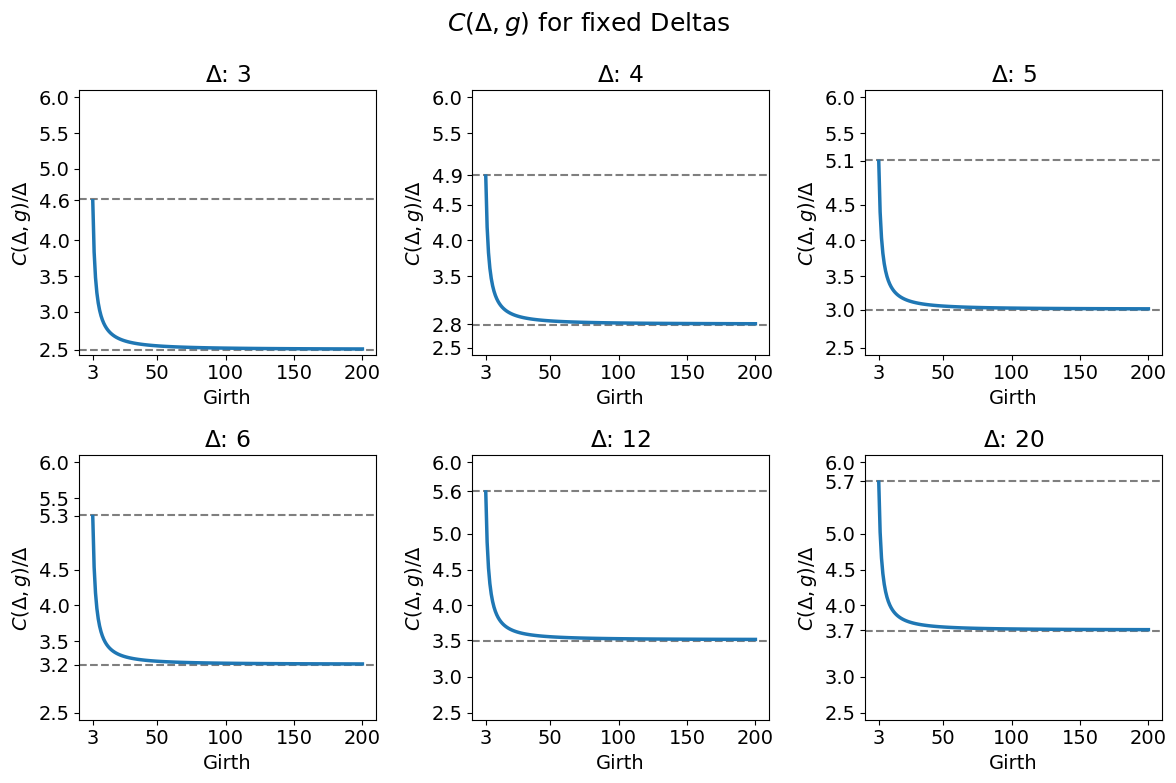}
\end{center}
\begin{center}
Figure 2. \footnotesize{ The quantity $C(\D,g)/\D$ as a function of $g\ge 3$ for  several fixed values of $\D$}
\end{center}
\end{figure}

\begin{figure}[!ht]
\begin{center}
\includegraphics[width=9cm,height=6cm]{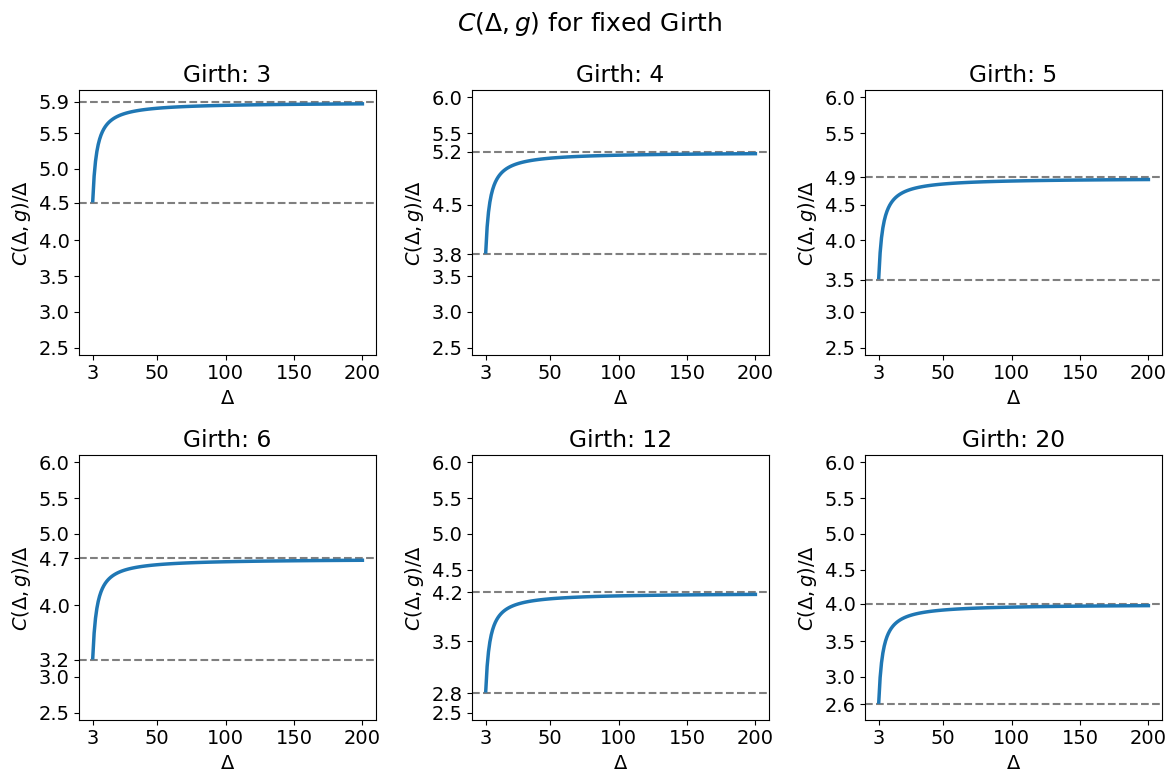}
\end{center}
\begin{center}
Figure 3. \footnotesize{The quantity $C(\D,g)/\D$ as a function of $\D\ge 3$ for  several fixed values of $g$.}
\end{center}
\end{figure}
\newpage

%
%
%
%
%
%
%
%

\numsec=2\numfor=1

\section{Proof of Theorem \ref{teowhit}}

\subsection{An alternative representation for $\Xi_\GI(q)$}

\\The present subsection, based in \cite{FJP}, aims to rewrite function $\Xi_\GI(q)$, defined in (\ref{XiG}) and (\ref{actrq}) in a
more convenient manner in terms of a sum over forests.  In order to do that, let us first establish  some notations  and recall some definitions regarding subgraphs, trees and forests of a graph $\GI$.

\\Given $n\in \N$, we set $[n]=\{1,2,\dots,n\}$ and given a finite set $X$, let $|X|$ denote its cardinality.
Along this paper, $\GI=(\VU, \EE)$ represents a fixed simple graph with vertex set $\VU$ and  edge set $\EE$.
A graph $\tilde\GI=(\tilde\VU,\tilde\EE)$ is a subgraph of $\GI=(\VU,\EE)$ if $\tilde \VU \subseteq \VU$ and $\tilde\EE \subseteq \EE$. A subgraph $\tilde\GI$
of $\GI$ is called spanning if $\tilde\VU=\VU$. A graph $\GI=(\VU, \EE)$
is {\it connected}
if for any pair $B, C$ of  subsets of $\VU$ such that
$B\cup C =\VU$ and $B\cap C =\emptyset$, there is at least an edge  $e\in \EE$ such
that $e\cap B\neq\emptyset$ and $e\cap C\neq\emptyset$.
A connected component of $\GI$  is a maximal connected subgraph of $\GI$.
A tree is a connected graph with no cycles.
A subset $R\subset \VU$ is said to be connected if $\GI|_{R}$ is connected; we denote by $\mathrm{C}_\VU$ the set of all connected subsets of $\GI$ with cardinality greater than one, i.e.
$ \mathrm{C}_\VU=\{R\subseteq \VU:\; \GI|_R \; \mbox{is connected and} \; |R|\ge 2 \}$.

\\We denote by $\mathrm{T}_\GI$ the set of (not necessarily spanning) trees which are subgraphs of $\GI$.
  An element $\t\in \mathrm{T}_\GI$ will be regarded as a collection of edges and we concur that the empty set is an element of $\mathrm{T}_\GI$. An element $\t\in \mathrm{T}_\GI$ distinct from the empty set is said to be non trivial. The number of edges of a non trivial tree $\t$ will be indicated by $|\t|$, and the set of vertices of a non trivial $\t$ will be denoted by $V_\t$. A forest in $\GI$ is a subgraph of $\GI$ whose connected components are trees.  We denote by $\mathfrak{F}_\GI$ the set of
all forests in $\GI$ and we agree that the empty set belongs to $\mathfrak{F}_\GI$.
An element $F\in \mathfrak{F}_\GI$ distinct from the empty set (i.e. a non trivial forest) is thus, for some $k\in \N$,  a collection $\{\t_1, \dots ,\t_k\}$ of non trivial trees such that
$V_{\t_i}\cap V_{\t_j}=\0$ for all pairs $\{i,j\}\subset [k]$. A non trivial forest $F$ of $\GI$  will
be  treated as a collection of edges. We denote by
$|F|$ the number of edges in $F$ and by $V_F$ the subset formed by all  vertices of $\GI$ belonging to the non trivial connected components of $F$. Namely, if $F$ is the collection
$\{\t_1, \dots ,\t_k\}$, then $V_F=\cup_{i=1}^kV_{\t_i}$ and $|F| = \sum_{i=1}^{k}{|\t_i|}$.

\\Given $R\in \mathrm{C}_\VU$, recall that $\mathcal{C}_{\GI|_R}$ represents the set of all connected spanning subgraphs of $\GI|_R$, and let us denote by
$\mathcal{T}_{\GI|_R}\subset \mathcal{C}_{\GI|_R}$ the set formed by  all  spanning trees of $\GI|_R$.

\\We  now rewrite the sum in the r.h.s. of \reff{actrq} using a beautiful and long known identity   discovered by O. Penrose in 1967 \cite{pen67}.
This identity  is based on  the existence of a map from  $\mathcal{T}_{\GI|_R}$ to $\mathcal{C}_{\GI|_R}$,
 called   {\it partition scheme} (see e.g. \cite{sok01}, \cite{scosok05} and \cite{JPS}).
\vv
\\\underline{\it Partition scheme in $\GI$.}
\\Given $R\in \mathrm{C}_\VU$, a partition scheme in $\GI|_R$ is a map ${ \mathrm{m}_{R}}:\mathcal{T}_{\GI|_R}\to\mathcal{C}_{\GI|_R}$
 such that for all $\t\in \mathcal{T}_{\GI|_R}$, it holds that
$\t\subset {\mathrm{m}_R}(\t)$ and
 $\mathcal{C}_{\GI|_R}=\biguplus_{\tau\in \mathcal{T}_{\GI|_R}}[\tau, \mathrm{m}_R(\tau)]$, where $\biguplus$ stands for disjoint union
 and $[\tau, \mathrm{m}_R(\tau)]=\{g\in \mathcal{C}_{\GI|_R}: \tau\subseteq g\subseteq {\mathrm{m}_R(\tau)}\}$.
 A partition scheme in $\GI$ is a family ${\mathrm{\mathbf{m}}}=\{\mathrm{m}_R\}_{R\in \mathrm{C}_\VU}$ such that
 for any $R\subset \mathrm{C}_\VU$, $\mathrm{m}_R$ is a partition scheme   on ${\cal C}_{\GI|_R}$. For any non trivial $\t$ subtree of $\GI$
 such that $V_\t=R$, we will denote shortly
$\mathrm{ \mathbf{m}}(\t)=\mathrm{m}_{R}(\t)$.

\\Given  a  partition scheme $\mi$  in $\GI$ and given $R\in \mathrm{C}_\VU$,
%
the Penrose identity  is the following equality
\be\label{ide}
\sum_{g\in\mathcal{C}_{\GI|_R}}(-1)^{|E_g|}= 
(-1)^{|R|-1} \sum_{\t\in \mathcal{T}_{\GI|_R}\atop \mi(\t)=\t}1.
\ee
See e.g. \cite{FP} for the proof.
Using  \reff{ide} we can rewrite
the activity $z(R,q)$ of a polymer $R$, given in \reff{actrq}, as follows
\be\label{actsche}
z(R,q)=\sum_{\t\in {\mathcal T}_{\GI|_R}\atop \mi(\t)=\t}\left(-{1\over q}\right)^{|\t|},
\ee
Inserting \reff{actsche} into \reff{XiG} we get that
$\Xi_\GI(q)$ can be rewritten
as
\be\label{whigen}
\Xi_\GI(q)=\sum_{F\in \mathfrak{F}^\mi_\GI} \left(-{1\over q}\right)^{|F|},
\ee
where the symbol  $\mathfrak{F}^{\mi}_\GI$ denotes the set of all forests in $\GI$ such that any non trivial tree $\t$ of the forest obeys
$\mi(\t)=\t$ and again  in the sum in the r.h.s. of \reff{whigen}   the empty set contributes with  the factor ``1".

\\In conclusion, recalling \reff{zqx}, we get that the chromatic polynomial $P_\GI(q)$ of a graph $\GI=(\VU,\EE)$ can be written as
\be\label{pgm}
P_\GI(q)= q^{|\VU|}\sum_{F\in \mathfrak{F}^\mi_\GI} \left(-{1\over q}\right)^{|F|}
\ee
where $\mi$ is any partition scheme in $\GI$.

\\It is remarkable to observe that,
although the set $ \mathfrak{F}^\mi$ depends strongly on  the partition scheme $\mi$ used,
 the identity \reff{pgm} holds independently of $\mi$.
Various partition schemes are available (see e.g. \cite{scosok05}, Sec. 2.2 and references therein), in particular, in applications on statistical mechanics and combinatorics,
as far as we know, two partition schemes have been used: the first one is the so called {\it Penrose partition scheme},  originally presented in \cite{pen67}; the second one  is  the {\it minimal-tree partition scheme} (described  e.g. in the proof of Lemma 2.2. in \cite{scosok05}). The Penrose  partition scheme  was used in \cite{FP} to improve the Dobrushin criterion for the convergence of the cluster expansion on the abstract polymer  and in \cite{BFPS} to improve the Lov\'asz local lemma.
The minimal tree partition scheme was used in \cite{PY} to improve the lower bound of the convergence radius of the Mayer Series on a gas of
classical particles in the continuum interacting via a stable and regular pair potential.

\\We will use  here the representation \reff{whigen} of the function $\X_\GI(q)$ with the partition scheme $\mi$ set to be equal the above mentioned minimal tree partition scheme, which we define here below.
\vv
\\\underline{\it Minimal tree partition scheme in $\GI$}.
Choose a total order $\succ$ in the set  $\EE$. For any $R\in \mathrm{C}_\VU$,  let $\mu_R: \mathcal{T}_{\GI|_R}\to\mathcal{C}_{\GI|_R}$ be the map that associates to each $\tau\in \mathcal{T}_{\GI|_R}$ the graph $\mu_R( \tau )\in \mathcal{C}_{\GI|_R}$
whose edge set is $E_{\mu_R(\tau)}= \tau \cup E_\tau^{*}$, where $E_\tau^{*}$ is composed by all edges $\{x,y\}\subset \EE|_R\setminus \t$
such that $\{x,y\} \succ \{z,u\} $ for every edge $\{z,u\} \in \tau$ belonging to the path from $x$ to $y$ through $\tau$.
It is easy to show that  $\mu_R$  is a partition scheme in $\mathcal{C}_{\GI|_R}$  (see e.g.  \cite{PY},\cite{PYob}, \cite{scosok05}), and thus $\bm \mu=\{\mu_R\}_{R\in \mathrm{C}_\VU}$ is a partition scheme in $\GI$.
Observe that
 for any tree $\t\subset \GI|_R$  the condition $\bm\m(\t)=\t$ simply means that for any $\{x,y\}\in \EE|_R\setminus \t$,
the path $p^\t_{xy}$ in $\t$ from $x$ to $y$ is such that $\{x,y\}$ is never the largest edge (in the chosen ordering $\succ$ in $\EE$) in the set
$\{\{x,y\}\cup \{e\}_{e\in p^\t_{xy}}\}$.

\\In graph theory a \emph{broken circuit} in the graph $\GI$  is  a set of edges of $\GI$ obtained by removing from some cycle in $\GI$  its largest edge (in the chosen order).
Therefore, the condition $\bm\m(\t)=\t$ in \reff{actsche}, with  $\bm\m$ being the minimal tree partition scheme, can be rephrased by saying that $\t$ does not contain broken circuits, that is to say, $\t$ is broken-circuit free (BCF).

\\In conclusion, denoting  by  $\mathfrak{F}^{\bm\m}_\GI$ the set of all forests in $\GI$  whose  non trivial trees are  broken-circuit free,  we get, by \reff{whigen}, that
\be \label{whitmu}
P_\GI(q)= q^{|\VU|} \sum_{F\in \mathfrak{F}^{\bm\m}_\GI} \left(-{1\over q}\right)^{|F|}.
\ee
The rewriting \reff{whitmu} of the chromatic polynomial $P_\GI(q)$ is known as
the Whitney's Broken Circuit Theorem \cite{W}. The discussion above provides therefore a ``statistical-mechanics" proof of this theorem.

\\From now on  we set $\mathcal{F}_\GI\equiv \mathfrak{F}^{\bm \mu}_\GI$, so that  $\mathcal{F}_\GI$  will denote the set of  all {BCF forests} of $\GI$.  Moreover $\mathcal{T}_\GI\subset \mathrm{T}_\GI$
will denote the set of  all { BCF trees} in $\GI$. Thus, a forest $F\in \mathcal{F}_\GI$  is either the empty set or,
for some  $k\in \N$, a set of BCF trees
$\{\t_1,\dots,\t_k\}$  such that $V_{\t_i}\cap V_{\t_j}=\0$ for all $\{i,j\}\subset [k]$.
\\With these notations, recalling  identity \reff{whitmu} the partition function $\Xi_\GI(q)$ can be rewritten  as
\be\label{xigf}
\Xi_\GI(q)=\sum_{F\in \mathcal{F}_\GI}\left(-{1\over q}\right)^{|F|}
\ee
We  emphasise once again that the empty forest is included in the sum of the r.h.s. of \reff{xigf} and it contributes with a factor 1.

\subsection{The zero-free region for $\Xi_\GI(q)$. Preliminary steps}

\\Let us thus investigate the zero-free region for $\Xi_\GI(q)$ starting from its representation given in \reff{xigf}.
 To achieve this objective we will follow the strategy based on induction developed
in ref. \cite{JPR}.

\\We start by  establishing some further notations.
 Given the graph $\GI=(\VU,\EE)$ and
given  $S \subseteq \VU$, let $\mathcal{F}_{\GI,S}$ be the set of all BCF forests of $\GI$  such that every tree  of $F$ contains at least one vertex of $S$ (again we agree that the empty set is an element of $\mathcal{F}_{\GI,S}$).
We write $\mathcal{F}^\bullet_{\GI,S}$ for those $F \in \mathcal{F}_{\GI,S}$ such that each tree  of $F$ contains exactly one vertex of $S$ (again, the empty set is included in  $\mathcal{F}^\bullet_{\GI,S}$). Clearly we have  that $\mathcal{F}^\bullet_{\GI,S} \subseteq \mathcal{F}_{\GI,S}\subseteq \mathcal{F}_{\GI}$.

\begin{definition}\label{def:functionF}
    For any $z\in \C$  and any graph $\GI$ we define the function
\be\label{FGz}
{F}_\GI(z) = \sum_{F \in \mathcal{F}_{\GI}}z^{|F|}.
\ee
\end{definition}
\vv

\\With an abuse of notation we will  also write ${F}_\VU(z)$ in place of $F_\GI(z)$, when we need to emphasize the dependency  on the vertex set. Moreover, for any $U \subseteq \VU$ and $S\subset U$ we will write shortly $\mathcal{F}_{U}$,
 $\mathcal{F}_{U,S}$  $\mathcal{F}^\bullet_{U,S}$ instead of $\mathcal{F}_{\GI|_U}$,  $\mathcal{F}_{\GI|_U,S}$ , $\mathcal{F}^\bullet_{\GI|_U,S}$  respectively, and, given  $\{u,w\}\subset U$ we denote by $\mathcal{T}_{U,v}$ and $\mathcal{T}_{U,v,w}$ the set of all BCF tree in $\GI|_U$ containing $v$ and containing $u$ and $v$ respectively. Note that $\mathcal{T}_{U,v}$ contains the empty set while
  $\mathcal{T}_{U,v,w}$  does not.

\\Let us list some  properties of function ${F}_\VU(z)$ which will be used  ahead. First,  given any  vertex $u\in \VU$, note that
\begin{equation}
\label{eq:recur}
{F}_\VU(z) = {F}_{\VU-u}(z) + \sum_{\substack{\t \in \mathcal{T}_{\VU,u} \\ |\t| \geq 1 }}z^{|\t|}{F}_{\VU-u -V_\t}(z)\, .
\end{equation}
Moreover, given $U\subseteq \VU$,  for any non empty set $S\subset U$, we have
\[
{F}_U(z)= \sum_{F\in \mathcal{F}_{U, S}}z^{|F|}{{F}_{U\setminus (S\cup V_F)}}(z).
\]
whence
\be\label{idun}
\sum_{F\in \mathcal{F}_{U,S}}z^{|F|}{{F}_{U\setminus (S\cup V_F)}(z)\over {F}_U(z)}=1, ~~~~\mbox{for all $U\subseteq \VU$ and $S\subset U$}.
\ee
Observe finally that by \reff{xigf} we have that
\be\label{xif}
{F}_\GI(-1/q)=\Xi_\GI(q),
\ee
\vv
\\According to identity \reff{xif},  in order to prove Theorem \ref{teowhit} we just need to find a zero-free region for the
function $F_\GI(z)$ in the complex plane. We will follow closely  the induction strategy  developed by Jensse, Patel and Regts in \cite{JPR}.
To this aim, in the following subsection we will extend the definition \reff{rodel2}-\reff{zdga} to the cases $\D=1$, $\D=2$ and $g=+\infty$ and
state  two lemmas needed to implement the induction.

\subsubsection{Extensions of definitions \reff{rodel2}-\reff{zdga} to the cases $\D=1$, $\D=2$ and $g=+\infty$ }
Given a graph $\GI=(\VU,\EE)$,
our goal is to  show that $F_\GI(z)\neq 0$
when $z$ belongs to a  disk in the complex plane centered at $z=0$ whose
radius depends
on the maximum degree $\D$  and the girth  $g$ of $\GI$ and the proof of this fact will be done by induction on the cardinality $|\VU|$ of the vertex set $\VU$ of $\GI$.  For the induction to work,
we
need  to extend the definitions \reff{rodel2}-\reff{zdga} to the cases $\D=1$, $\D=2$ and $g=+\infty$. This is necessary because
we will prove the statement above for a given graph $\GI=(\VU,\EE)$ with maximum degree $\D$ and girth $g$ having assumed by induction that the statement is true for all graphs   $\tilde\GI=(\tilde \VU,\tilde \EE)$  such that $|\tilde \VU|<|\VU|$ (and thus, in particular, for all subgraphs  of $\GI$). Now,
any subgraph $\tilde\GI$ of $\GI$ has necessarily
maximum degree $\tilde \D\le \D$ and girth $\tilde g\ge g$, so that  we could have $\tilde \D=1$, $\tilde\D=2$ and $\tilde g=+\infty$ for some
$\tilde \GI\subset \GI$.

\\Let us thus define
\begin{equation}\label{RD12}
\r_\D= \begin{cases}
+\infty &{\rm if}~ \D=2\\
+\infty &{\rm if}~ \D=1
\end{cases}
~~~~~~{\rm and}~~~~~~
R_\D= \begin{cases}
1 &{\rm if}~ \D=2\\
+\infty &{\rm if}~ \D=1
\end{cases}
\end{equation}
and for any $\D\ge 1$ let
 \be\label{finf}
f^d_\D(b)=0 ~~~~\mbox{if $d=+\infty$}.
 \ee

 \\With the extension \reff{RD12} and \reff{finf}, we can see that the function  $x_\D(b)$  given in \reff{hdb} is, for any $\D\ge 1$, an increasing function of $b\in[1,\r_\D)$, varying from 0 to $R_\D$ and that definitions \reff{hdb}-\reff{zdga} are now valid  also for the cases  $\D=1$, $\D=2$ and $g=+\infty$.

 \\In the special case in which $\GI$ is a graph with maximum degree $\D=1$ (and thus with its girth is $g$ necessarily infinite) we get that $z_1^\infty (a)=a$ for all $a\in [0,1)$ and  we define conventionally $z_1^\infty(1):=0$ so that
 \be\label{zia}
 z_1^\infty(a)= \begin{cases} a &{\rm if~} a\in [0,1)\\
 0 &{\rm if~} a=1
 \end{cases}
 \ee

 \\For the case $\D\ge 2$ and $g$ either finite or infinite we have the following Lemma.

\begin{lemma}\label{lemma:RemarkB}
For any fixed $a\in [0,1]$, for any $\D\ge 2$ and for any $3\le g<+\infty$, the quantity $b_\D^g(a)$ is  the unique solution of the equation
$K^g_\D(a,b)= a$ as $b$ varies in the interval $[1,\r_\D)$.  Moreover, if $g = \infty$, then  $b_\D^g(a)$ and  $z_\D^g(a)$ are correctly  defined by \reff{badg2} and \reff{zdga}.
\end{lemma}
\\{\bf Proof}.  When $g$ is finite and $\D\ge  3$ the first claim of Lemma \ref{lemma:RemarkB},
coincides with the thesis of Lemma \ref{lemma:RemarkA} which have been already proved. When $\D=2$ and $3\le g<\infty$
 the quantity $K^g_\D(a,b)$ given in \reff{acca}
 is again, for any fixed  $a\in [0,1]$, a continuous increasing function in the interval   $b\in[1,\r_\D)$  such that
  $K^g_\D(a,1)=0$ and $\lim_{b\to \r_\D}K^g_\D(a,b)=+\infty$. Therefore $b_\D^g(a)$ continues to be  the unique solution of the equation
$K^g_\D(a,b)= a$ as $b$ varies in the interval $(1,\r_\D)$.

\\On the other hand if $g=+\infty$ and $\D\ge 2$, we have, according to \reff{RD12} and \reff{finf},
\begin{equation*}
K^{\infty}_\D(a,b)= \Big[1+{(1-a)}x_\D(b)\Big]^{\Delta} - 1
\end{equation*}
which, for any fixed $a\in [0,1]$, is still non  decreasing in the interval $b\in [1,\r_\D)$ but now
$$ \lim_{b\to\r_\D}K^{\infty}_\D(a,b)=(1+(1-a)R_\D)^\D-1.
$$
Therefore, if $a_\D$ is the solution of the equation
\be\label{eqad}
(1+(1-y)R_\D)^\D-1=y,
\ee
then  the equation $K^\infty_\D(a,b)=a$ has solution only if $a<a_\D$ while for   $a\in [a_\D,1]$ we have that $K^\infty_\D(a,b)<a$ in the whole interval $b\in [1,\r_\D)$. Hence,
according to  definitions
\reff{badg2}, \reff{zdga},
we get
\be\label{zdinfa}
z_\D^{\infty}(a)=\begin{cases} (a+1)^{1\over \D}-1 & {\rm if}~ a\in [0,a_\D)\\
(1-a)R_\D & {\rm if}~ a\in [ a_\D,1]
\end{cases}
\ee
~~~~~~~~~~~~~~~~~~~~~~~~~~~~~~~~~~~~~~~~~~~~~~~~~~~~~~~~~~~~~~~~~~~~~~~~~~~~~~~~~~~~~~~~~~~~~~~~~~~~~~~~~~~~~~~~~~~~~~~~~~~$\Box$

\\Letting $\Ni=\N\cup{+\infty}$ and $\Ni_{\ge 3}=\Ni\setminus\{1,2\}$, observe that by \reff{zdinfa} and \ref{zia}) we have
\[
z_\D^g(a)|_{a=0}=z_\D^g(a)|_{a=1}=0~~~~~~ \mbox{for all $\D\in \N$ and all $g\in \Ni_{\ge 3}$}
\]

\begin{lemma}\label{incre}
For any fixed $a\in [0,1]$ the quantity $z_\D^g(a)$ given in \reff{zdga} is a decreasing sequence of $\D\in  \N$   and an increasing sequence of
$g\in \Ni_{\ge 3}$. Namely, for any $\D\ge 1$ and any $g\in \Ni\setminus\{1,2\}$, it holds that
\[
z_{\D}^g(a)\ge z_{\D+1}^g(a)
\]
\[
 z_{\D}^{g}(a)\le z_\D^{g+1}g(a)
\]
\end{lemma}
{\bf Proof}.
\underline{\it Case ${g=+\infty}$}. Let us thus consider initially the case in which  $\GI$ has no cycles and thus has girth $g=+\infty$.
We need to prove that, for any fixed $a\in [0,1]$, we have $z_{\D}^\infty(a)\ge z_{\D+1}^\infty(a)$ for any $\D\ge1$. Indeed, first observe that according to \reff{zia} and \reff{zdinfa} we have that
$z_1^\infty(a)\ge z_\D^\infty(a)$ for all $\D\ge 2$.  Then note that the sequence
$\{a_\D\}_{\D\ge 2}$ (where recall that $a_\D$ is the solution of the equation \reff{eqad}) is decreasing and such that $a_2={1\over 2}(5-\sqrt{13})$ and
$$
\lim_{\D\to \infty}a_\D=eW(e^{-1+e/2})-1\approx 0,295741
$$
where $W$ is the product log function.
Now, according to \reff{zdinfa},   $z_\D^{\infty}(a)$ is, for any $\D\ge 2$, a continuous function in the closed interval $a\in [0,1]$,  starting at the value 0 when $a=0$, increasing in the interval $[0,a_\D]$, reaching   its maximum at $a=a_\D$ and then decreasing
in the interval $(a_\D,1]$ reaching the value zero at $a=1$.
The fact that, for any $\D\ge 2$, we have that  $a_\D\ge a_{\D+1}$  and that $(a+1)^{1\over \D}-1\ge (a+1)^{1\over \D+1}-1$ for all $a\ge 0$,
implies that  and for any $\D\ge 2$ we have that the inequality $z_{\D}^\infty(a)\ge z_{\D+1}^\infty(a)$ holds in the whole interval $a\in [0,1]$.

\vv
\\\underline{\it Case ${g<+\infty}$ and $\D\ge 2$}. In this case we have  that $z_\D^g(a)=(1-a)x_\D(b_\D^g(a))$
where $b_\D^g(a)$ is the solution of the equation $K_\D^g(a,b)=a$ in the interval $b\in[1,\r_\D)$. Recalling that
 $z_\D^g(0)=z_\D^g(1)=0$ for  all $\D\in \N$ and all $g\in \Ni_{\ge 3}$, we need to prove that
 $z_\D^g(a)\ge z_{\D+1}^g(a)$ when $a\in (0,1)$.
 Since $x_\D(b)$  is increasing for $b\in [1,\r_\D)$, and decreasing with $\D$ for any fixed $b>1$, it is sufficient to prove that
$b_\D^g(a)$ defined in \reff{badg2}  is a decreasing sequence of $\D$ and an increasing sequence of $g$ for any $a\in (0,1)$ and for any $\D\ge 2$.
To see this,  recall that  $K_\D^g(a, b)$ is a continuous monotonic increasing  function of $b$ in the interval $[1,\r_\D)$ such that  $K_\D^g(a, 1)=0$ and $\lim_{b\to \r_\D}K_\D^g(a, b)=+\infty$.
Moreover we have that:

\\{\it (i)} $\r_{\D+1}<\r_\D$;

\\{\it (ii)}  $K_{\D+1}^g(a, b)>K_{\D}^g(a, b)$
for any fixed $a\in (0,1)$ and for all $b\in (0,\r_{\D+1})$.

\\Items {\it (i)} and {\it (ii)} imply that $b_{\D+1}^g(a)\le b_\D^g(a)$ for  all $a\in (0,1)$ and all $\D\ge 2$ and thus  $z_{\D+1}^g(a)\le z_\D^g(a)$ for  all $a\in (0,1)$ and all $\D\ge 2$.

\\Finally, it is  immediate  to see that $b_{\D}^{g}(a)\le  b_\D^{g+1}(a)$ for any $a\in (0,1)$ since for any fixed $a\in (0,1)$ and all
$b\in (0,\r_\D)$ it holds  that  $K_{\D}^g(a, b)>K_{\D}^{g+1}(a, b)$ (because $(\D-1)(1-b^{-1\over \D})<1$ for any
$b\in (0,\r_\D)$ when $\D\ge 2$).

 \\It remains to prove that, for all finite $g\ge 3$ and for all $\D\ge 2$, we have that $z_\D^g(a)\le z_\D^\infty(a)$. This is because, by the discussion above regarding the case $g=\infty$,  by the fact that $K^g_\D(a,b)> K^\infty_\D(a,b)$  for all $g\in \N\setminus\{1,2\}$ and by the fact that $\lim_{b\to \r_\D}K^g_\D(a,b)=+\infty$, we have that $b^{g}_\D(a)\le b^{\infty}_{\D}(a)$ for any finite $g\ge 3$ and thus $z^{g}_\D(a)\le z^{\infty}_{\D}(a)$. \underline{}

 ~~~~~~~~~~~~~~~~~~~~~~~~~~~~~~~~~~~~~~~~~~~~~~~~~~~~~~~~~~~~~~~~~~~~~~~~~~~~~~~~~~~~~~~~~~~~~~~~~~~~~~~~~~~~~~~~~~~~~~~~~~~~~~~~~~~~~~~~~~~~~~~~$\Box$

\\In the next subsection we will state  one  more lemma  (Lemma \ref{le2}) concerning generation functions of trees and we conclude by enunciating  a proposition (Proposition \ref{pro1})
about the zero free region of $F_\GI(z)$ from which Theorem \ref{teowhit} follows straightforwardly.  The proof of Lemma \ref{le2}, which is lengthly, will be given in the appendix, while Section \ref{secpro} will be devoted to the proof of Proposition \ref{pro1}.

\subsubsection{A lemma about generation functions of trees and a proposition about the zero free region of $F_\GI(z)$}
Recall that  $\mathrm{T}_\GI$ denotes  the set of (not necessarily spanning) trees which are subgraphs of $\GI=(\VU,\EE)$ and that an element $\t\in \mathrm{T}_\GI$ is regarded as  a set of edges with $|\t|$   denoting the number of edges of $\t$. We also remind that empty set is an element of $\mathrm{T}_\GI$. Given  a vertex $v\in \VU$, we let $\mathrm{T}_{\GI,v}$ be
the subset of $\mathrm{T}_\GI$  formed by those  trees   in $\mathrm{T}_\GI$ whose vertex set contains $v$. We also agree that
 $\0\in\mathrm{T}_{\GI,v}$.
We define for $x\ge 0$
\be\label{tgv}
T_{\GI,v}(x):=\sum_{\substack{\t\in \mathrm{T}_{\GI,v}}} x^{|\t|}.
\ee
Note that the sum in the r.h.s. of \reff{tgv} is a polynomial in $x$ starting with    the term 1 which occurs when $\t=\0$.

\\Analogously,  given $\{v,w\}\subset \VU$ and
denoting by  $\mathrm{T}_{\GI,v,w}$
the subset of $\mathrm{T}_\GI$  formed by those trees whose vertex set  contains $v,w$,
we define, for $x\ge 0$ the function
\be\label{tgv12}
T_{\GI,v,w}(x)= \sum_{\t\in \mathrm{T}_{\GI,v,w}} x^{|\t|}
\ee
Observe that  $\mathrm{T}_{\GI,v,w}$ is now a polynomial in $x$ with starting the term $c_dx^d$  where $d$ is the length of the shortest
path between $v$ and $w$ in $\GI$ and $c_d$ is the number of paths of length $d$  between $u$ and $v$ in $\GI$.

\begin{lemma}\label{le2} Given a graph  $\GI=(\VU,\EE)$  with maximum  degree $\D\ge 2$ and girth $g$,
let $d$ be the length of the shortest path between two given vertices $v$ and $w$ in $\GI$. Let $T_{\GI,v}(x)$ and $T_{\GI,v,w}(x)$ be the functions defined in \reff{tgv} and \reff{tgv12} respectively, then, for any $b\in (1, \r_\D)$ with $\r_\D$ defined in \reff{rodel2}
and for any non-negative $x\le x_\D(b)$ with $x_\D(b)$ defined in \reff{hdb},
we have that
\be\label{ineqv}
T_{\GI,v}(x)\le b
\ee
and
\be\label{ineqvw}
T_{\GI,v,w}(x)\le   bf^d_\D(b)
\ee
where  $d$  is the number of edges of the shortest path between $v$ and $w$ in $\GI$ and $f^d_\D(b)$ is the function defined in \reff{fdelta2}.
\end{lemma}
As anticipated above, the proof of Lemma \ref{le2} is given  in  Appendix A.

\\We conclude this section stating   a proposition (whose proof requires Lemmas \reff{incre} and \reff{le2})
from which Theorem \ref{teowhit} follows straightforwardly.
For this purpose, we introduce an auxiliary function of a complex variable.
\vv
\begin{definition}
    \\Given a vertex $u\in \VU$  and given
$z\in \C$ such that  ${F}_{\VU-u}(z)\neq 0$, we define
\be\label{rug}
R^u_\GI(z)= \frac{{F}_{\VU}(z)}{{F}_{\VU-u}(z)} - 1
\ee
and if $\GI=(\VU,\EE)$ is such that  $\VU=\{u\}$, i.e. $\VU-u=\0$,  we take $R^u_\GI(z)=0$.
\end{definition}

\vskip.2cm
\\Recalling identity  \reff{eq:recur},   the function $R^u_\GI(z)$ can also be written,   as follows
\begin{equation}
\label{RuG}
R^u_\GI(z)= \sum_{\substack{\t \in \mathcal{T}_{\VU,u} \\ |\t| \geq 1 }}z^{|\t|}
\frac{{F}_{\VU-u-V_\t}(z)}{{F}_{\VU-u}(z)}\, ,
\ee
where we recall that $\mathcal{T}_{\VU,u}$ is the set of all BCF trees of $\GI$ containing $u$.

\vv

\begin{proposition}\label{pro1}
Let $\GI=(\VU,\EE)$ be a graph with maximum degree  $\D\ge 1$ and girth $g\in \Ni_{\ge 3}$.
Given   $a \in [0,1]$, let  $z^g_\D(a)$  be the quantity  defined  \reff{zdga}.
Then,
for  any  $z\in \C$ satisfying $|z| \leq z_\D^g(a)$,
    we have that
\be\label{claim1}
{F}_\GI(z)\neq 0.
\ee
and
\be\label{rlea}
|R^u_\GI(z)|\le a
\ee
where $z^g_\D(a)$, ${F}_\GI(z)$ and $R^u_\GI(z)$ are the functions defined in \reff{zdga}, \reff{FGz} and \reff{rug} respectively.
\end{proposition}

\vv

\subsection{Proof of Proposition \ref{pro1}}\label{secpro}
We will prove Proposition \ref{pro1} by induction on the cardinality $|\VU|\ge 2$, following closely the work done in \cite{JPR}.
\vv
\\- When $|\VU|=1$ then ${F}_\GI(z)=1$ and $R^u_\GI(z)=0$ for all $\ \in \mathbb{C}$.
\vv
\\- If $\GI=(\VU,\EE)$ is a graph such that $|\VU|=2$,  then either $\D=0$ or $\D=1$. Suppose that $\D=0$, then again $F_\GI(z)=1$ and $R^u_\GI(z)=0$ for all $\ \in \mathbb{C}$. Now, if $|\VU|=2$ and $\D=1$ (and consequently $g=+\infty$) we have that  ${F}_\GI(z)=1+z$,
and $R^u_\GI(z)=z$ and $z_1^\infty(a)$ is given by \reff{zia} and therefore claims \reff{claim1} and \reff{rlea} hold.

\\In conclusion,  Proposition \ref{pro1} is true for all graphs
 $\GI=(\VU,\EE)$ with $|\VU|\le 2$.
 \vv
 \\Furthermore  observe that  any graph $\GI=(\VU,\EE)$  with maximum degree $\D=1$  and $\EE\neq\0$ (and consequently $g=\infty$) is a collection of isolated vertices and a collection of  non intersecting isolated edges  so that  $F_\GI(z)=(1+z)^k$ where
 $k\ge 0$ is the number of  isolated edges of $\GI$. Hence again we get that $F_\GI(z)\neq 0$ if $|z|\le z_1^\infty(a)$ for all $a\in [0,1]$.
 To bound $|R^u_\GI(z)|$ observe that $F_\GI(z)/F_{\GI-u}(z)$ is either equal to 1 if $u$ is an isolated vertex of $\GI$ or it is equal
 to $1+z$ if $u$ belong to an isolated edge of $\GI$. Thus, we have  either $R^u_\GI(z)=0$ or $R^u_\GI(z)=z$, then in any case
 $|R^u_\GI(z)|\le a$ for any $|z|\le z_1^\infty(a)$ with $a\in [0,1]$.
\vv
\\- Now let $\GI=(\VU,\EE)$ be a graph  with $|\VU|=n > 2$, maximum degree  $\D$  and girth  $g$  (by the last observation here above
 we can assume $\D\ge 2$) and assume that Proposition \ref{pro1} is true for all graphs  $\widetilde\GI=(\widetilde\VU,\widetilde\EE)$ with $|\widetilde\VU|<|\VU|$.
 \vv

\\Given $u\in \VU$,   let us set shortly $\GI'=(\VU-u, \EE|_{\VU-u})$ and  $\VU'=\VU-u$. Then by the induction hypothesis, given $a\in (0,1)$, we have that $F_{\VU'}(z)\neq 0$ for $|z|\le z_{\D'}^{g'}(a)$ where
 $\D'$ and $g'$ are the maximum degree and the girth of $\GI'$ respectively. Since $\D'\le \D$ and $g'\ge g$, we get,
 by Lemma \ref{incre},  that  $ z_{\D'}^{g'}(a)\ge z_{\D}^{g}(a)$ and thus $F_{\VU-u}(z)\neq 0$ for $|z|\le z_{\D}^{g}(a)$.
In light of this, to conclude the proof by induction of Proposition \ref{pro1},  it is sufficient to prove    that, given $u \in \VU$ and  any $a\in (0,1)$  
\be\label{aux}
|R^u_\GI(z)|\le a~~~~~~~\mbox{for all $z\in \C$ such that  $|z|\leq z_\D^g(a)$}.
\ee
Indeed, by the induction hypothesis and Lemma \ref{incre} we have that $|{F}_{\VU-u}(z)|>0$ for $|z|\leq z_\D^g(a)$. Therefore
 \reff{aux}  immediately implies
$$
|R^u_\GI(z)|\le a~~\Longrightarrow~~\Big|\frac{{F}_\VU(z)}{{F}_{\VU-u}(z)} - 1 \Big|\le
a~~\Longrightarrow~~\Big|\frac{{F}_\VU(z)}{{F}_{\VU-u}(z)}\Big|\ge 1-a > 0
~~\Longrightarrow~~{F}_\VU(z)\neq 0.
$$

\\To show  \reff{aux}, we need to use  repetitively that, for any  $A\subset \VU'$,
\be\label{indu1}
\left| \frac{{F}_{\VU'-A}(z)}{{F}_{\VU'}(z)} \right|\leq (1-a)^{-|A|}~~~~~\mbox{ for any $|z|\le z_\D^g(a)$}.
\ee
Indeed, bound \reff{indu1} is a consequence of the induction hypothesis. To see this,
let  $A = \{v_1, \ldots, v_k\}  \subseteq \VU'$ and for $i\in [k]$, let $A_i = \{v_1, \ldots, v_i \}$ and $A_0 = \emptyset$.
Let $\D_i$ and $g_i$ be the maximum degree and the girth of the graph $\GI|_{\VU'-A_i}$.
Then by the induction hypothesis, for any $i\in [k]$ and for any  $|z|\le  z_{\D_i}^{g_i}(a)$, we have that ${F}_{\VU'-A_{i}}(z)\neq 0$ and
\be\label{indu2}
\left| \frac{{F}_{\VU'-A_{i}}(z)}{{F}_{\VU'-A_{i+1}}(z)} - 1 \right|\le a~~\Longrightarrow~~ 
\left|\frac{{F}_{\VU'-A_{i+1}}(z)}{{F}_{\VU'-A_{i}}(z)} \right|\le {1\over 1-a}~.
\ee

\\Now observe that $\D_i$ is decreasing with $i$  (i.e. $\D_{i+1}\le \D_{i}$ for all $i\in [k]$) and $g_i$ is increasing with $i$
(i.e. $g_{i+1}\ge g_{i}$).
Therefore,  by  Lemma \ref{incre}, the sequence $ z^{g_i}_{\D_i}(a)$ is increasing with $i$ and in particular
$ z^{g_i}_{\D_i}(a) \ge z_\D^g(a)$ for all $i\in [k]$. Hence,  if   $|z|\le z_\D^g(a)$ we get
that
 \reff{indu2}
surely holds true for all $i\in [k]$ .
So we get
\[
\left| \frac{{F}_{\VU'-A}(z)}{{F}_{\VU'}(z)} \right|
= \prod_{i=1}^{k} \left| \frac{{F}_{\VU'-A_{i}}(z)}{{F}_{\VU'-A_{i-1}}(z)} \right|
\leq (1-a)^{-|A|},~~ ~~~~ \forall z\in  \C: ~ |z|\le z_\D^g(a).
\]

\\To conclude the proof of Proposition \ref{pro1} we now just need to  prove the inequality  \reff{aux}.
Let us denote by $\G_{\GI}(u)=\{v\in \VU: \{u,v\}\in \EE\}$  the (punctured) vertex neighbor of $u$ in $\GI$ and by
$E_{\GI}(u)=\{e\in \EE: \{u\}\subset e\}$ the set of edges of $\GI$ incident to $u$.
We  assume that the edges in $E_u(\GI)$ are  the largest in the ordering of $\EE$ chosen to construct the minimal tree partition scheme $\bm \mu$. We stress that
such choice of the ordering $\prec$ in $\EE$ is without loss of generality  because for any $U\subset \VU$ the function $F_U(z)$ does not depend on the order established in $\EE$ (actually not even on the partition scheme chosen!). For later convenience, we  also establish an order in $\G_u(\GI)$ in such a way that,
for any $\{v,w\}\subset \G_u(\GI)$ such that  $\{u,v\}\prec \{u,w\}$ we have that $v<w$.

\\We now have all ingredients to proof the inequality  \reff{aux}, i.e., to prove that  $|R^u_\GI(z)|\le a$ for any $\GI=(\VU,\EE)$ with maximum degree $\D\ge 2$ and girth $g\in \bar\N\setminus\{1,2\}$ when $|\VU|> 2$ under the hypothesis that $|z|\le z_\D^g(a)$. We start by  the expression of  $R^u_\GI(z)$ given in \reff{RuG} and rewrite it  as follows:
\begin{equation}
\label{Rugiz}
\begin{aligned}
R^u_\GI(z)
& =\sum_{S \subseteq \G_\GI(u) \atop S \not= \emptyset} z^{|S|}
\sum_{F \in \mathcal{F}^\bullet_{\VU',S} \atop  F \cup uS \text{ is BCF}}z^{|F|}
\frac{F_{\VU\setminus(S\cup V_F)}(z)}{{F}_{\VU'}(z)}\\
\end{aligned}
\ee

\\where  $uS$ is a short notation for   $\{\{u,s\}: s \in S\}$ and we recall that $\mathcal{T}_{\GI,u}$ is the set of all BCF trees contained  in $\GI$ and containing $u$.
 The  equality    \reff{Rugiz} holds because, due to our choice of the order in $ \EE$ (i.e. if  $e\in E_\GI(u)$ then $e\succ f$, for all $f\in \EE\setminus E_\GI(u)$),  any tree $\t\in \mathcal{T}_{\GI,u}$  is
the union of the edges of $\t$ incident to $u$, say $uS$ for some  $S\subset \G_\GI(u)$, with a BCF forest $F\in \mathcal{F}^\bullet_{\VU',S} $ under the condition that this union is BCF. From \reff{Rugiz} we now have
\[
\begin{aligned}
R^u_\GI(z)& =\sum_{\substack{S \subseteq \G_\GI(u) \\ S \not= \emptyset}} z^{|S|}
\sum_{\substack{{F} \in \mathcal{F}^\bullet_{\VU',S}  }}z^{|F|} \frac{F_{\VU\setminus(S\cup V_F)}(z)}{{F}_{\VU'}(z)}
-
\sum_{\substack{S \subseteq \G_\GI(u) \\ S \not= \emptyset}} z^{|S|}
\sum_{\substack{{F} \in \mathcal{F}^\bullet_{\VU',S} \\ F \cup uS \text{ is not  BCF} }}z^{|F|} \frac{{F}_{\VU' - S - V_F}(z)}{{F}_{\VU'}(z)}\\
& = \sum_{\substack{S \subseteq \G_\GI(u) \\ S \not= \emptyset}} z^{|S|}
\sum_{\substack{{F} \in \mathcal{F}_{\VU',S}  }}z^{|F|} \frac{F_{\VU\setminus(S\cup V_F)}(z)}{{F}_{\VU'}(z)}
-
\sum_{\substack{S \subseteq \G_\GI(u) \\ S \not= \emptyset}} z^{|S|}
\sum_{\substack{{F} \in \mathcal{F}^2_{\VU',S} }}z^{|F|} \frac{F_{\VU\setminus(S\cup V_F)}(z)}{{F}_{\VU'}(z)}~~-
\\
&~~~~-
\sum_{\substack{S \subseteq \G_\GI(u) \\ S \not= \emptyset}} z^{|S|}
\sum_{\substack{{F} \in \mathcal{F}^\bullet_{\VU',S} \\ F \cup uS \text{ is not  BCF} }}z^{|F|}
\frac{F_{\VU\setminus(S\cup V_F)}(z)}{{F}_{\VU'}(z)}\\
\end{aligned}
\]
where in the second line  $\mathcal{F}^2_{\VU',S}$ is the set whose elements are those forests $F \in \mathcal{F}_{\VU', S}$ where some component of $F$ contains (at least) two vertices of $S$.

\\Recalling  now \eqref{idun}, i.e.,
$$
\sum_{\substack{{F} \in \mathcal{F}_{\VU',S}  }}z^{|F|} \frac{F_{\VU\setminus(S\cup V_F)}(z)}{{F}_{\VU'}(z)}=1
$$
we get
\be
\label{Rugiz3}
\begin{aligned}
R^u_\GI(z) & = \sum_{\substack{S \subseteq \G_\GI(u) \\ S \not= \emptyset}} z^{|S|}
-
\sum_{\substack{S \subseteq \G_\GI(u) \\ |S| \ge 2}} z^{|S|}
\sum_{\substack{{F} \in \mathcal{F}^2_{\VU',S}  }}z^{|F|} \frac{F_{\VU\setminus(S\cup V_F)}(z)}{{F}_{\VU'}(z)}
\\
&-
\sum_{\substack{S \subseteq \G_\GI(u) \\ S \not= \emptyset}} z^{|S|}
\sum_{\substack{{F} \in \mathcal{F}^\bullet_{\VU',S} \\ F \cup uS \text{ is not  BCF} }}z^{|F|} \frac{F_{\VU\setminus(S\cup V_F)}(z)}{{F}_{\VU'}(z)}\\
\end{aligned}
\ee

\\Now observe that, by  the induction hypothesis and bound~\reff{indu1}, we have
\be\label{modK}
\left|z^{|F|}\frac{F_{\VU\setminus(S\cup V_F)}(z)}{{F}_{\VU'}(z)}\right|= |z|^{|F|} \left| \frac{F_{\VU\setminus(S\cup V_F)}(z)}{F_{\VU'}(z)} \right|\le |z|^{|F|} (1-a)^{-|V_F \cup S|}\le
\left( \frac{|z|}{1-a}\right)^{|F|}  (1-a)^{-|S|},
\ee
\vv
\\where the last inequality in \reff{modK} follows because
$|V_F \cup S| \leq |F| + |S|$ for all $F \in \mathcal{F}_{\VU',S} $.
Hence, setting  $\frac{|z|}{1-a}=x$, from \reff{Rugiz3} we get
\[
\begin{aligned}
|R^u_\GI(z)| & \le \sum_{S\subseteq \G_\GI(u) \atop S \not= \emptyset} [(1-a)x]^{|S|}
+
\sum_{S \subseteq \G_\GI(u) \atop |S| \ge 2}x^{|S|}
\sum_{F \in \mathcal{F}^2_{\VU',S}}x^{|F|}
\\
&+
\sum_{\substack{S \subseteq \G_\GI(u) \\ S \not= \emptyset}} x^{|S|}
\sum_{\substack{{F} \in \mathcal H}^\bullet_{\VU',S}}x^{|F|}  \\
\end{aligned}
\]
where we have set shortly
${\mathcal H}^\bullet_{\VU',S}=\{F\in  \mathcal{F}^\bullet_{\VU',S}:  F \cup uS \text{ is not  BCF} \}$.
Then we can bound
\begin{equation}\label{eq:RP1P2}
\begin{aligned}
|R^u_\GI(z)| & \le [1+(1-a)x]^\D-1 +P_1(x)+P_2(x)
\\
\end{aligned}
\end{equation}
where
\be\label{p1x}
P_1(x)=
\sum_{\substack{S \subseteq \G_\GI(u) \\ |S|\ge 2}}x^{|S|}
\sum_{{F} \in \mathcal{F}^2_{\VU',S} }x^{|F|}
\ee
and
\be\label{p2x}
P_2(x)=
\sum_{S \subseteq \G_\GI(u) \atop |S| \ge 1}x^{|S|}
\sum_{{F} \in  {\mathcal H}_{\VU',S} }x^{|F|}~.\\
\ee

\\We now  bound   $P_1(x)$ given in \reff{p1x}.
By the considerations above, we know that any forest $F\in \mathcal{F}^2_{\VU',S} $ is the union of  a
tree $\t\in \mathcal{T}_{\VU',v,w}$  where $\{v,w\}\subset S$ with a forest $F\in \mathcal{F}_{\VU'\setminus V_\t, S\setminus V_\t}$. Therefore,
we can bound
\be\label{p1xbound}
\begin{aligned}
\sum_{{F} \in \mathcal{F}^2_{\VU',S} }x^{|F|}&\le \sum_{\{v,w\}\subseteq S}
\sum_{\t\in  \mathcal{T}_{\VU',v,w}}x^{|\t|}
\sum_{F\in \mathcal{F}_{\VU'\setminus V_\t, S\setminus V_\t}} x^{|F|}\\
&\le  \sum_{\{v,w\}\subseteq S}
\sum\limits_{\t\in  \mathcal{T}_{\VU',v,w}}x^{|\t|}\sum_{F\in \mathcal{F}_{\VU', S\setminus \{v,w\}}}x^{|F|}\\
&\le  \sum_{\{v,w\}\subseteq S}
\sum\limits_{\t\in  \mathcal{T}_{\VU',v,w}}x^{|\t|}
\prod_{v'\in  S\setminus \{v,w\}}
\sum_{\s\in \mathcal{T}_{\VU', v'}}x^{|\s|}\\
&\le  \sum_{\{v,w\}\subseteq S}
\sum\limits_{\t\in  \mathrm{T}_{\VU',v,w}}x^{|\t|}
\prod_{v'\in  S\setminus \{v,w\}}
\sum_{\s\in \mathrm{T}_{\VU', v'}}x^{|\s|}\\
&\leq \sum_{\{v,w\}\subseteq S}  T_{\GI', v, w}(x) \prod_{v' \in S \setminus\{v,w\}}T_{\GI',v'}(x)\\
&\leq \sum_{\{v,w\}\subseteq S}  T_{\GI', v, w}(x) \prod_{v' \in S \setminus \{v\}} T_{\GI',v'}(x)\\
\end{aligned}
\ee
where the last inequality follows from the fact that $T_{\GI',v'}(x)\ge 1$.
Observe that $T_{\GI', v, w}(x)=0$ if there is no path in $\GI'$ connecting $v$ to $w$.  If there are paths in $\GI'$ connecting $v$ to $w$  then the shortest path from $v$ to $w$ in $\GI'$ plus the edges $\{u,v\}$, $\{u,w\}$ form a cycle of
length $dist(v,w)+2\ge g$ in $\GI$ so that $dist(v,w)\ge g-2$.
Now, by
hypothesis  we have that  $|z|\le z_\D^g(a)=(1-a)x_\D^g(a)$  so that $x\le x_\D^g(a)$. Therefore,  by  inequalities \reff{ineqvw} and \reff{ineqv} in Lemma \ref{le2}, we have that, for any $x\le x_\D(b)$
\be\label{Ftree}
T_{\GI', s, t}(x)\le \max_{\{u,v\}\subset \VU\atop dist(u,v)= g-2} T_{\GI, u, v}(x)\le bf^{g-2}(b).
\ee
and
\be\label{Ftree2}
T_{\GI',s'}(x)  \le  \max_{u\in \VU} T_{\GI, u}(x)\le b
\ee
so that
\be\label{Ftree3}
\sum_{{F} \in \mathcal{F}^2_{\VU',S} }x^{|F|}\le\sum_{\{v,w\}\subseteq S}f_\D^{g-2}(b)b^{|S|}=f_\D^{g-2}(b)b^{|S|}\sum_{\{v,w\}\subseteq S}1~.
\ee
Therefore
we can bound, for $x\le x_\D(b)$,
$$
\begin{aligned}
P_1(x)
&\le  {f_\D^{g-2}(b)}\sum_{S \subseteq \G_\GI(u) \atop|S|\ge 2}[bx_\D(b)]^{|S|}\sum_{\{v,w\}\subseteq S}1\\
&= {f_\D^{g-2}(b)}\sum_{S \subseteq \G_\GI(u) \atop|S|\ge 1}[bx_\D(b)]^{|S|}\sum_{\{v,w\}\subseteq S}1
\end{aligned}
$$
where in the last line if $|S|=1$ then $\sum_{\{v,w\}\subseteq S}1=0$.
\vv
\\We now  bound   $P_2(x)$ given in \reff{p2x}.  To this end observe that
if  $F\in {\mathcal H}^\bullet_{\VU',S}$ then,
$F$ must necessarily contain a tree
$\t\in\mathcal{T}_{\VU',v,w}$  where now $v\in S$ and $w\in \G_\GI(u)\setminus S$ such that $\{u,v\}\prec \{u,w\}$ (for the proof of this fact see Proposition 2.4 of \cite{JPR}).

\\So, recalling that we have ordered $\G_\GI(u)$ such that,
for $\{v,w\}\subset \G_\GI(u)$, $\{u,v\}\prec\{u,w\}\Rightarrow v<w$ and proceeding analogously to what was done along steps  \reff{p1xbound}-\reff{Ftree3},  we can bound
\vv
\[
\begin{aligned}
\sum_{{F} \in \mathcal{H}^\bullet_{\VU ',S}}x^{|F|}
&\le \sum_{v\in S}\sum_{ w\in \G_\GI(u)\setminus S\atop w>v}
\sum_{\t\in  \mathcal{T}_{\VU',v,w}}x^{|\t|}
\sum_{F\in \mathcal{F}_{\VU'\setminus V_\t, S\setminus V_\t}} x^{|F|}\\
&\le \sum_{v\in S}\sum_{ w\in \G_\GI(u)\setminus S\atop w>v}
\sum\limits_{\t\in  \mathcal{T}_{\VU',v,w}}x^{|\t|}\sum_{F\in \mathcal{F}_{\VU', S\setminus \{v\}}}x^{|F|}\\
&\le  \sum_{v\in S}\sum_{ w\in \G_\GI(u)\setminus S\atop w>v}
\sum\limits_{\t\in  \mathcal{T}_{\VU',v,w}}x^{|\t|}
\prod_{v'\in  S\setminus \{v\}}
\sum_{\t\in \mathcal{T}_{\VU', v'}}x^{|\t|}\\
&\leq \sum_{v\in S}\sum_{ w\in \G_\GI(u)\setminus S\atop w>v} T_{\GI', v, w}(x) \prod_{v' \in S \setminus \{v\}}T_{\GI',v'}(x)\\
&\le f_\D^{g-2}(b)b^{|S|}\sum_{v\in S}\sum_{ w\in \G_\GI(u)\setminus S\atop w>v}1
\end{aligned}
\]
where the last inequality holds for all $x\le x_\D(b)$.
Hence, for $x\le x_\D(b)$
\vv
$$
\begin{aligned}
P_1(x)+P_2(x)&\le {f_\D^{g-2}(b)}\sum_{S \subseteq \G_\GI(u) \atop |S|\ge 1}[bx_\D(b)]^{|S|}
\Big[\sum_{\{v,w\}\subset S}1+\sum_{v\in S}\sum_{ w\in \G_\GI(u)\setminus S\atop w>v}1\Big]\\
&=  {f_\D^{g-2}(b)}\sum_{S \subseteq \G_\GI(u) \atop |S|\ge 1}[bx_\D(b)]^{|S|}
\sum_{v\in S}\sum_{w\in \G_\GI(u)\atop w>v} 1
\end{aligned}
$$
\vv
\\Now let  us set $|\G_\GI(u)|=n$  (where of course $n\le \D$), and  set also $\G_\GI(u)=\{v_1,v_2,\dots, v_n\}$ such that
$\{u,v_1\}\prec\{u,v_2\}\prec\dots\{u,v_n\}$ and thus
$v_1<v_2<\dots<v_n$ according the order chosen in $\G_\GI(u)$,  then we have

$$
\begin{aligned}
\sum_{S \subseteq \G_\GI(u) \atop |S|\ge 1}[bx_\D(b)]^{|S|} \sum_{v\in S}\sum_{w\in \G_\GI(u)\atop w>v}1&=
\sum_{k=1}^n[bx_\D(b)]^{k}\sum_{J \subseteq [n] \atop |J|=k} \sum_{i\in J}\sum_{ j\in [n]\atop j>i}1\\
&=\sum_{k=1}^n[bx_\D(b)]^{k}\sum_{J \subseteq [n] \atop |J|=k} \sum_{i\in J}(n-i)\\
& =\sum_{k=1}^n[bx_\D(b)]^{k}\sum_{i\in [n]}(n-i)\sum_{J \subseteq [n]: i\in J \atop |J|=k}1 \\
&= bx_\D(b)\sum_{k=1}^{n}[bx_\D(b)]^{k-1}{n\choose 2}{n-1\choose k-1}\\
&={n\choose 2}bx_\D(b)\sum_{r=0}^{n-1}[bx_\D(b)]^{r}{n-1\choose r} \\
&={n\choose 2}bx_\D(b)(1+ bx_\D(b))^{n-1}\\
\end{aligned}
$$
\vv
\\i.e., for any $\le x_\D(b)$, recalling that $n\le \D$,  we have
\\
$$
P_1(x)+P_2(x)\le {f_\D^{g-2}(b)}{\D\choose 2}bx_\D(b)(1+ bx_\D(b))^{\D-1}.
$$

\\In conclusion, recalling (\ref{eq:RP1P2}), for any $z\in \C$ such that  ${|z|\over 1-a}=x\le x_\D(b)$ we obtain
\begin{equation}\label{eq:Rfinal}
\begin{aligned}
|R^u_\GI(z)| & \le  [1+(1-a)x_\D(b]^\D-1+{f_\D^{g-2}(b)}{\D\choose 2}bx_\D(b)(1+ bx_\D(b))^{\D-1}\\
&= K_\D^g(a,b)
\end{aligned}
\end{equation}
and choosing $b=b_\D^g(a)$, for $|z|\le (1-a)x_\D(b_\D^g(a))= z_\D^g(a)$ we get the bound
$$
|R^u_\GI(z)|\le a
$$
which concludes the induction and therefore the proof of Proposition \ref{pro1}.

~~~~~~~~~~~~~~~~~~~~~~~~~~~~~~~~~~~~~~~~~~~~~~~~~~~~~~~~~~~~~~~~~~~~~~~~~~~~~~~~~~~~~~~~~~~~~~~~~~~~~~~~~~~~~~~~~~~~~~~~~~~~~~~~~~~~~~~~~~~$\Box$

\subsection{Conclusion of the proof of Theorem \ref{teowhit}}
Theorem \ref{teowhit} now  follows straightforwardly from Proposition \ref{pro1}.
Indeed, Proposition \ref{pro1} states that  for  any   $a\in [0,1]$ and any $z\in \C$ satisfying $|z| \leq z_\D^g(a)$ with $\D\in \N$ and $g\in \Ni_{\ge 3}$, we have that ${F}_\GI(z)\neq 0$. Now,  according to  Lemma \ref{lemma:RemarkA} in section \ref{sec3}, $z_\D^g(a)$ is a continuous function of $a$ in the interval $[0,1]$
such that $z_\D^g(0)=z_\D^g(1)=0$ and therefore admits a maximum. Choosing for $a$ the value at which $z_\D^g(a)$
is maximum we get that  ${F}_\GI(z)\neq 0$ as soon as $|z| \leq \max_{a\in [0,1]}z_\D^g(a)$.
Finally,
since $\Xi_\GI(q)= {F}_\GI(-1/q)$, this implies  that for any   $q\in \C$ satisfying
$|q| \geq \big[\max_{a\in (0,1)}[z_\D^g(a)]\big]^{-1}\equiv C(\D,g)$,
the partition function  $\Xi_\GI(q)\neq 0$, and thus, $P_\GI(q)$ is
free of zeros when  $|q| \geq C(\D,g)$.

\section{ Proof of Theorem \ref{coro22}}
Let us first observe that, for any $\D\ge 3$ and any $b>1$, we have
$$
x_\D(b)
\ge {\ln b\over \D b}
$$
and, since $\inf_{\D\ge1}\r_\D = e$, with $\r_\D$ given in \reff{rodel2}, by Lemma \ref{le2} we can bound, for all $b\in[1, e)$,
$$
T_{\GI,v}\left({\ln b\over \D b}\right)\le b
~~~~~~~~~~~~~~~\mbox{and}~~~~~~~~~~
T_{\GI,v_1,v_2}\left({\ln b\over \D b}\right)\le  bf^d_\D(b),
$$
 Moreover, observe that when $b\ge 1$
$$
1-b^{-{1\over \D}}\le {\ln b\over \D}
$$
and therefore, when $g$ is finite, by \reff{fdelta2}
$$
\begin{aligned}
f^{g-2}_\D(b)
& \le{{1\over (\D-1)}~{[(\D-1) {\ln b\over \D}]^{g-2}\over  [1-(\D-1) {\ln b\over \D}]}}\\
&\le {{1\over (\D-1)}~{[(\D-1) {\ln b\over \D}]^{g-2}\over  [1-{\ln b}]}}\\
&\le {1\over \D -1} {(\ln b)^{g-2}\over  1-{\ln b}}.
\end{aligned}
$$
Recalling the bound for $|R^u_\GI(z)|$ given in the first line of \reff{eq:Rfinal}, for  any $\D\ge 3$ and any $g\ge $ finite, taking $|z|\le {(1-a)\ln b\over \D b}$ with $b\in [1,e)$, we get$$
\begin{aligned}
|R^u_\GI(z)|
&\le \left[1+ { (1-a)\ln b\over \D b}\right]^\D-1
+ {f^{g-2}_\D(b)} \binom{\D}{2}\left({\ln b\over \D }\right) \left(1+{\ln b\over \D }\right)^{\D-1}\\
&\le \exp\left\{{(1-a)\ln b\over b}\right\} -1
+ {f^{g-2}_\D(b)}  \binom{\D}{2}\left({\ln b\over \D }\right)\left(1+{\ln b\over \D }\right)^{\D}\\
&\le \exp\left\{{(1-a)\ln b\over b}\right\} -1
+ {{1\over \D-1} {(\ln b)^{g-2}\over  1-{\ln b}}}  \binom{\D}{2}\left({\ln b\over \D }\right)\left(1+{\ln b\over \D }\right)^{\D}\\
&\le \exp\left\{{(1-a)\ln b\over b}\right\} -1
+  {b(\ln b)^{g-1}\over  2(1-{\ln b})} \\
&\doteq K^g_\infty(a,b)
\end{aligned}
$$
where in the penultimate line we have used that $(1+{\ln b\over \D })^{\D}\le b$ for all $\D\ge 1$.
Note that $K^g_\infty(a,b)$ is  for any $a\in (0,1)$ a monotonic increasing function of $b$ in the interval $[1,e)$ such that
$K^g_\infty(a,1)=0$ and $\lim_{b\to e} K^g_\infty(a,b)=+\infty$.
Thus for $b\in (1, e)$ and $a\in [0,1]$,  letting $b_\infty^g(a)$ be the (unique) solution of the equation $K^g_\infty(a,b)=a$, we get that
$P_\GI(q)$
is free of zeros if
$$|q|\ge K_g\D
$$
where
\[
K_g=\left[\max_{a\in [0,1]}\left({(1-a)\ln b_\infty^g(a)\over b_\infty^g(a) }\right)\right]^{-1}
\]

\\Now,
for any fixed $a\in [0,1]$, we have that  $K^g_\infty(a,b)\ge K^g_\D(a,b)$ as $b\in [1, e)$. Therefore we have that $b_\infty^g(a)\le b_\D^g(a)$
for any $\D\ge 3$ and thus ${\ln b_\infty^g(a)\over \D b_\infty^g(a)}\le x_\D(b_\D^g(a))$ for any $a\in [0,1]$.
Therefore, for all $a\in (0,1)$, we  have that
\be\label{aia}
\max_{a\in [0,1]}{(1-a)\ln b_\infty^g(a)\over \D b_\infty^g(a)}\le \max_{a\in [0,1]}z_\D^g(a).
\ee
and  \reff{aia} implies that
$$
{C(\D, g)\over \D}\le K_g ~~~~{\rm for~ all}~~ \D\ge 3.
$$

\\To see that $C(\D,g)/\D$ is increasing with $\D\ge 3$ for any finite $g\ge 3$, according to definition \reff{CDg} we  need to show that
$$
\D\big[ \max_{a\in [0,1]} z_\D^g(a)\big]\ge \max_{a\in [0,1]}(\D+1)\big[z_{\D+1}^g(a)\big].
$$
To this end, we recall that, according to the proof of Lemma \ref{incre},  $b_\D^g(a)$ is decreasing with $\D$ for any $a\in [0,1]$.
Recalling Lemma \ref{lemma:RemarkA}, $b_\D^g(a)$  is continuous in the interval $a\in [0,1]$.
Moreover   $x_\D(b)$ given in \reff{hdb} is increasing with $b$ for any fixed $\D\ge 3$ and decreasing with
$\D$ for any fixed $b\in [1,\r_\D)$. These two facts
imply  that the sequence $\{\D x_\D(b_\D^g(a))\}_{\D\ge 3}$ is decreasing sequence of continuous functions  of  $a\in [0,1]$.  Therefore,
setting $f_\D^g(a)=\D z_\D^g(a)= (1-a)\D x_\D(b_\D^g(a))$, we have that   $\{f_\D^g(a)\}_{\D\ge 3}$ is a decreasing sequence of continuous of $a\in [0,1]$. Observe also that by  construction $f_\D^g(0)=0$ and also $f_\D^g(0)=1$, so that $f_\D^g(a)$ admits maximum in $[0,1]$, so that
$\{\max_{a\in [0,1]}f_\D^g(a)\}_{\D\ge 3}$ is a decreasing sequence of $\D$ and thus $C(\D,g)/\D= \{(\max_{a\in [0,1]}f_\D^g(a))^{-1}\}_{\D\ge 3}$  is an increasing
sequence of $\D$.

\section*{Acknowledgements}
It is a pleasure to  thank Benedetto Scoppola for discussions and  comments. P.M.S.F. and  E.J.  were supported by
FAPEMIG (Funda\c{c}\~ao de Amparo \`a Pesquisa do Estado de Minas Gerais).

\section*{Competing Interests Statement}
The authors declare none.

%

\section*{Data Availability Statement}
No data was used for the research described in the article.

\section*{Funding statement}
This work was supported by  Funda\c{c}\~ao de Amparo \`a Pesquisa do Estado de Minas Gerais (FAPEMIG)
under research grant  34026 (Bolsa de Desenvolvimento em Ci\^encia, Tecnologia e Inova\c{c}\~ao, N\'\i vel $\setminus$I Inc.). The funder had no role in study design, data collection and analysis, decision to publish, or preparation of the manuscript.

\renewcommand{\theequation}{A.\arabic{equation}}
\setcounter{equation}{0}

\section*{Appendix A: Proof of Lemma \ref{le2}}

\underline{\bf Case $\D\ge 3$}.
\vskip.1cm

\\Let $T(\D)$  be the $\D$-regular infinite tree. Given
$x\ge 0$,   consider the following function
\be\label{tx}
T(x)=\sum_{n\ge 1} t_{n}(\D) x^{n}
\ee
where $t_{n}(\D)$ is the number of  subtrees of  $T(\D)$ with $n$ vertices and containing a fixed vertex. Then, if $\GI$ is a graph with maximum degree $\D$,  we have that
\be\label{tvt}
T_{\GI,v}(x)\le {T(x)\over x}
\ee
Let us also define $U(\D)$ the infinite rooted tree with $\D-1$ children in each vertex  and let $u_n(\D)$ be the set of $n$-vertices trees
contained in $U(\D)$ and rooted with the same root as $U(\D)$. Let
\be\label{ux}
u(x)=\sum_{n\ge 1} u_{n}(\D) x^{n}
\ee
Let us finally define $U^*(\D)$ the infinite rooted tree with $\D-1$ children in each vertex except for the root which have $\D-2$ children
and let $u^*_n(\D)$ be the set of $n$-vertices trees
contained in $U^*(\D)$ and rooted with the same root as $U^*(\D)$.
Let
\be\label{usx}
u^*(x)=\sum_{n\ge 1} u^*_{n}(\D) x^{n}
\ee
Then we have that
\be\label{Us}
u(x)=x(1+u(x))^{\D-1}
\ee
and therefore
\[
u^*(x)=x(1+u(x))^{\D-2}
\]
\[
T(x)=x(1+u(x))^{\D}
\]
Define:
\[
W(x)={u(x)\over x}= (1+u(x))^{\D-1}=\sum_{n\ge 1} u_{n}(\D) x^{n-1}
\]
\[
W^*(x)={u^*(x)\over x}= (1+u(x))^{\D-2}=\sum_{n\ge 1} u^*_{n}(\D) x^{n-1}
\]
\be\label{sx}
S(x)={T(x)\over x}= (1+u(x))^{\D}=\sum_{n\ge 1} t_{n}(\D) x^{n-1}
\ee

\\{\bf Proposition A.1}. {\it The series (\ref{tx}), (\ref{ux}), (\ref{usx}) converge for $x\in [0,R)$ where
\begin{equation*}
R={(\D-2)^{\D-2}\over (\D-1)^{\D-1}}
\end{equation*}
Moreover, for
$b\in (1, ({\D-1\over \D-2})^\D)$ and
\be\label{xb2}
x_\D(b)={b^{1\over \D}-1\over b^{\D-1\over \D}}= {b^{2\over \D}-b^{1\over \D}\over b}
\ee
we have that
\be\label{wws}
W^*(x_\D(b))= b^{1-{2\over \D}}~, ~~~~~~W(x_\D(b))= b^{1-{1\over \D}},~~~~~~~ S(x_\D(b))= b
\ee
}
\\{\bf Proof}.
Observe  that $x_\D(b)$ is an increasing function as $b$ varies  in the interval
$b\in (1, ({\D-1\over \D-2})^\D)$  such that $x_\D(1)=0$ and $x_\D(({\D-1\over \D-2})^\D)=R$.
The relation (\ref{Us}) define $u$ as an implicit function of $x$ and implies that
$$
u=f^{-1}(x)
$$
where
\begin{equation*}
f(u)\;=\;{u \over (1+u)^{\D-1}}
\end{equation*}
For $u\ge 0$,  the function $f$ starts at $0$ and it is  strictly increasing in the interval
$[0,(\D-2)^{-1}]$, at the end of which it attains its maximum value $f((\D-2)^{-1})=R$ with
$$
R= {(\D-2)^{\D-2}\over (\D-1)^{\D-1}}
$$
Therefore, $f$ is a bijection from $[0,(\D-2)^{-1}]$ onto $[0,R]$ and so $u(x)$ satisfying (\ref{Us}) is defined
in the whole interval $x\in (0.R)$. By Lagrange inversion formula we have that
$$
u_n(\D) ~=~{1\over n!} {d^{n-1}\over du^{n-1}}(1+u)^{(\D-1)n}\Big|_{u=0}~={1\over n!} {[(\D-1)n]!\over [(\D-2)n+1]!}
$$
and thus the convergence radius of the series (\ref{ux})
is, as expected,
$$
R=\lim_{n\to\infty}{u_n(\D)\over u_{n+1}(\D)}= {(\D-2)^{\D-2}\over (\D-1)^{\D-1}}
$$
The series \reff{usx} and (\ref{tx}) are also converging in the same convergence radius because
$u^*_n(\D)\le u_n(\D)$ for all $n\in \mathbb{N}$ and $T(x)=u(x)+u^2(x)$.

\\To show \reff{wws}, observe that $S(x)=\frac{1}{x} T(x)=(1+u(x))^\D$ is an increasing function in the interval $(0,R)$ varying from $1$
to $({\D-1\over \D-2})^\D$ and therefore, for any $b\in (1,({\D-1\over \D-2})^\D)$ we have that
$$
S^{-1}(b)= {b^{1\over \D}-1\over b^{\D-1\over \D}}=x_\D(b)
$$ Indeed
$$
S(x)= b~~\Rightarrow~~[1+u(x)]^{\D}= b~~\Rightarrow~~u(x)= b^{1\over \D}-1 ~~\Rightarrow~~x= f( b^{1\over \D}-1)
~~\Rightarrow~~x= {b^{1\over \D}-1\over b^{\D-1\over \D}}
$$
Moreover, since $u(x_\D(b))=b^{1\over \D}-1$, we get
$$
W(x_\D(b))= (1+u(x_\D(b)))^{\D-1}= b^{\D-1\over \D}
$$
and
$$
W^*(x_\D(b))= (1+u(x_\D(b)))^{\D-2}= b^{\D-2\over \D}
$$

~~~~~~~~~~~~~~~~~~~~~~~~~~~~~~~~~~~~~~~~~~~~~~~~~~~~~~~~~~~~~~~~~~~~~~~~~~~~~~~~~~~~~~~~~~~~~~~~~~~~~~~~~~~~~~~~~~~~~~~~~~~~~~~~~~~~~~~~~~~~~~~~~~$\Box$
\vv
\\Recalling now \reff{tvt} and \reff{sx} and using Proposition  (A.1)  we get
$T_{\GI,v}(x_\D(b))\le S(x_\D(b))=b$
and thus  also
$$
T_{\GI,v}(x)\le S(x_\D(b))\le b, ~~~~~{\rm for~ all}~x\le x_\D(b)
$$
which proves inequality \reff{ineqv} when $\D\ge 3$.
To prove inequality \reff{ineqvw}, let $p(m)$ be the number of paths of length $m$ (i.e. with $m$ edges) with end points $v$ and $w$
and let $d$ be  distance between $v$ and $w$ (i.e. the number of edges of the shortest path from $u$ to $w$).

\\Consider first the case $d=1$.
Then, for any  $x\le x_\D(b)$, we get
$$
\begin{aligned}
T_{\GI,v,w}(x) & \le x(W(x))^2+ \sum_{m\ge 2}p(m)x^m(W^*(x))^{m+1} \\
& =x(W(x))^2+ \sum_{m\ge 2}(\D-1)^{m-1}x^m(W^*(x))^{m+1} \\
& = x(W(x))^2+ {W^*(x)\over \D-1} \sum_{m\ge 2}(\D-1)^{m}x^m(W^*(x))^{m}\\
& = x(W(x))^2+ {(\D-1)(W^*(x))^3 x^2\over 1- (\D-1)x(W^*(x))}\\
&\le {b(1-b^{-1\over \D})\over 1-(\D-1)(1-b^{-1\over \D})}
\end{aligned}
$$
where in the last line we have used that $x\le x_\D(b)$ and the first two of the three identities \reff{wws}.

\\Let us now consider  the case  $dist(v,w)=d>1$. Then for  $x\le x_\D(b)$ we get
\be\label{tolong}
\begin{aligned}
T_{\GI,v,w}(x)& \le x^d(W(x))^2(W^*(x))^{d-1} + \sum_{m\ge d+1}p(m)x^m(W^*(x))^{m+1}]\\
&\le x^d(W(x))^2(W^*(x)^{d-1} + {W^*(x)\over \D-1}\sum_{m\ge d+1}[(\D-1)x^mW^*(x)]^{m}]\\
& =  x^d(W(x))^2(W^*(x)^{d-1} + {W^*(x)\over \D-1}{[(\D-1)xW^*(x)]^{d+1}\over 1- (\D-1)xW^*(x)}\\
& \le b(1-b^{-1\over \D})^d + {{b^{1-{2\over \D}}}(\D-1)^d (1-b^{-1\over \D})^{d+1}\over 1- (\D-1)(1-b^{-1\over \D})}
\end{aligned}
\ee
where in the last inequality we have use the hypothesis  $x\le x_\D(b)$ and  identities \reff{wws}. Now, via a simple algebraic manipulation of the r.h.s. of the last line of \reff{tolong} we can bound
\be\label{tlo2}
\begin{aligned}
T_{\GI,v,w}(x) &\le  b(1-b^{-1\over \D})^d + {{b^{1-{2\over \D}}}(\D-1)^d (1-b^{-1\over \D})^{d+1}\over 1- (\D-1)(1-b^{-1\over \D})}\\
&= {b[(\D-1)(1-b^{-1\over \D})]^d\over 1- (\D-1)(1-b^{-1\over \D})}
\left[{1\over (\D-1)^d}+ (1-b^{-1\over \D})\Big(b^{-{2\over \D}}-{1\over (\D-1)^{d-1}} \Big)\right]\\
& \le {b[(\D-1)(1-b^{-1\over \D})]^d\over 1- (\D-1)(1-b^{-1\over \D})}
\left[{1\over (\D-1)^d}+ {1\over \D-1}\Big(1-{1\over (\D-1)^{d-1}} \Big)\right]\\
&={b\over \D-1}{[(\D-1)(1-b^{-1\over \D})]^d\over 1- (\D-1)(1-b^{-1\over \D})}
\end{aligned}
\ee
where in the penultimate line of \reff{tlo2} we have used that $1-b^{-1\over \D}\le {1\over \D-1}$ in the interval $b\in [1, \big({\D-1\over \D-2}\big)^{\!\D})$ and that $b^{-2\over \D}\le 1$, so that ${1\over (\D-1)^{d}}-(1-b^{-{1\over \D}}){1\over (\D-1)^{d-1}}\le 0$ and
$(1-b^{-1\over \D})b^{-{2\over \D}}\le {1\over \D-1}$. In conclusion we get that, for any $\{u,w\}\subset \VU$ and for any $d\ge 1$,
$$
T_{\GI,v,w}(x)\le {b\over \D-1}{[(\D-1)(1-b^{-1\over \D})]^d\over 1- (\D-1)(1-b^{-1\over \D})}=bf^d_\D(b)
$$
The proof of inequality \reff{ineqvw} is therefore concluded for the case $\D\ge 3$.
\vv

\\
\vskip.3cm
\\\underline{\bf CASE $\D=2$}
\vskip.1cm

\\When $\D=2$ we have  that
$$
u(x)= \sum_{n\ge 1} x^n={x\over 1-x}
$$
and thus
$$
T(x)=\sum_{n\ge 1}nx^n= x(1+u(x))^2={x\over (1-x)^2}
$$
Both series $u(x)$ and $T(x)$ has convergence radius $R=1$. Therefore
$$
S(x)= {1\over (1-x)^2}
$$
is a bijection from $[0,1)\to [1,+\infty)$. Let  $b\in [1,+\infty)$,
then $S(x) = b$ implies that $x= 1-b^{-1\over 2}$. I.e.
$$
x_{\D=2}(b)=1-b^{-1\over 2}= {b^{2\over 2}-b^{1\over 2}\over b}
$$
consistently with \reff{xb2}. Therefore when $\GI$ is a graph with maximum degree $\D=2$ we get for any $x\le x_\D(b)$,
$$
T_{\GI,v}(x)\le S(x_\D(b))=b
$$
and if $v$ and $w$ are vertices of $\GI$ at distance $d$ then
$$
T_{\GI,v,w}(x)\le {((x_\D(b))^d\over (1-(x_\D(b))^2}= b(1-b^{-{1\over 2}})^d\le  b^{3/2}(1-b^{-{1\over 2}})^d= b f^d_2(b).
$$
which concludes the proof of Lemma \ref{le2}.

\renewcommand{\theequation}{B.\arabic{equation}}
\setcounter{equation}{0}


%


\begin{thebibliography}{99}


\bibitem{BFP} R. Bissacot, R. Fern\'andez, A. Procacci: {\it On the convergence of cluster expansions for polymer gases}, {J.  Stat. Phys.}, {\bf 139}, 598-617  (2010).

\bibitem{BFPS} Bissacot, R.; Fern\'andez, R.; Procacci A.; Scoppola, B.: {\it An Improvement of the Lov\'asz Local Lemma via Cluster Expansion},
Comb. Probab.  Comput., {\bf 20}, n. 5, 709-719
(2011).




\bibitem{dob96} R. L. Dobrushin: {\it Estimates of semiinvariants
  for the {I}sing model at low temperatures}.  {Topics in
    Statistics and Theoretical Physics}, Amer. Math. Soc. Transl. (2),
  {\bf 177}, 59-81 (1996).


\bibitem{FP} R. Fernandez; A. Procacci: {\it Cluster expansions for
abstract polymer models.  New bounds from an old approach}, Comm. Math.
Phys., {\bf 274},  123-140 (2007).


\bibitem{FP2} R. Fern\'andez, A. Procacci: {\it Regions Without Complex Zeros for Chromatic Polynomials on Graphs with Bounded Degree},
Comb. Probab. and Comput., {\bf 17},  225-238 (2008).

\bibitem{FJP} P. M. S. Fialho, E. Juliano, A. Procacci: {\it A remark on the Whitney Broken Circuit Theorem}. arXiv:2407.04035,  https://arxiv.org/abs/2407.04035
(2024).


\bibitem{JPS} B. Jackson; A. Procacci; A. D. Sokal: {\it Complex zero-free regions at large $|q|$ for multivariate Tutte
polynomials (alias Potts-model partition functions) with
general complex edge weights}, J.  Comb. Theor. B, {\bf 103}, 21-45 (2013).


\bibitem{JPR} M. Jenssen, V. Patel, G. Regts: {\it Improved bounds for the zeros of the chromatic polynomial via Whitney's Broken Circuit Theorem}. J.  Comb. Theor. B,
 {\bf 169}, 233-252 (2024).

\bibitem{KE} D. Kim; I.G. Enting: {\it The limit of chromatic polynomials}, J.  Comb. Theory B,
{\bf 26}, Issue 3, 327-336 (1979).



 \bibitem{pen67} O. Penrose: Convergence of fugacity
 expansions for classical systems.  In {\it Statistical
 mechanics: foundations and applications}\/, A. Bak (ed.),
 Benjamin, New York (1967).

 \bibitem{PY}   A. Procacci and S. A. Yuhjtman: {\it  Convergence of Mayer and virial expansions
and the Penrose tree-graph identity}, Lett. Math. Phys., {\bf 107}, 31-46 (2017).

\bibitem{PYob}  A. Procacci and S. Yuhjtman: {\it The Mayer series and the Penrose tree gaph identity}, Oberwolfach Rep. 14, in: Miniworkshop: Cluster expansions: From Combinatorics to Analysis through Probability (R. Fernandez, S. Jansen, D. Tsagkarogiannis, eds.) (2017).




\bibitem{scosok05}  A. Scott, A. D. Sokal: {\it The repulsive lattice gas, the independent-set polynomial, and the Lov\'asz local lemma}. {J. Stat. Phys.}  {\bf 118}, n. 5-6, 1151--1261 (2005).

\bibitem{sok01} A. Sokal: {\it Bounds on the complex zeros of
(di)chromatic polynomials and Potts-model partition functions}.
{Combin. Probab. Comput.,} {\bf 10}, 41-77  (2001).

 \bibitem{W} H. Whitney (1932): {\it A logical expansion in mathematics}. { Bull. American Math. Soc.},  {\bf 38}, n. 8, 572-579.



\end{thebibliography}
\end{document}